\documentclass[11pt,onecolumn]{article}
\usepackage[left=2.5cm,top=2.4cm,right=2.5cm,bottom=2.4cm]{geometry}
\usepackage{graphicx}
\usepackage{amssymb,amsmath,amsthm,amsopn,mathrsfs,amsfonts}
\usepackage{times}
\usepackage{multirow,caption,enumerate}
\graphicspath{{./Figures/}}	

\newcounter{remarkcounter} 
\theoremstyle{definition}
\newtheorem{remark}[remarkcounter]{Remark}
\newtheorem{theorem}{Theorem}[section]
\newtheorem{corollary}[theorem]{Corollary}
\newtheorem{lemma}[theorem]{Lemma}

\DeclareMathOperator{\eps}{\varepsilon}
\newcommand{\s}{\mathfrak{s}}
\newcommand{\ts}{\bar{\s}}

\title{\vspace{-20pt}\Large{From Canards of Folded Singularities to \\Torus Canards in a Forced van der Pol Equation}}
\author{
John Burke\footnotemark[2]\,
{}\footnotemark[4], 
Mathieu Desroches\footnotemark[3],
Albert Granados\footnotemark[3], \\
Tasso J. Kaper\footnotemark[2],
Martin Krupa\footnotemark[3],
and Theodore Vo\footnotemark[2]
}

\begin{document}
\maketitle

\renewcommand{\thefootnote}{\fnsymbol{footnote}}
\footnotetext[2]{Department of Mathematics and Statistics, Boston University, 111 Cummington Mall, Boston, MA 02215, USA}  
\footnotetext[3]{INRIA Paris-Rocquencourt Research Centre, MYCENAE Project-Team, Domaine du
Voluceau, Rocquencourt BP105, 78153 Le Chesnay cedex, France}
\footnotetext[4]{Now at MSCI, Inc., 2100 Milvia Street, Berkeley, CA, 94704}
\renewcommand{\thefootnote}{\arabic{footnote}}
\let\thefootnote\relax\footnote{\textbf{Running Title}: From Canards of Folded Singularities to Torus Canards}
%=========================================================================================

%------------------------------------------------------------	
\vspace{-10pt}
\begin{abstract}        \label{sec:abstract}
\noindent
In this article, we study canard solutions of the forced van der Pol equation in the relaxation limit for low-, intermediate-, and high-frequency periodic forcing. A central numerical observation made herein, which motivated our study, is that there are two branches of canards in parameter space which extend across all positive forcing frequencies. In the low-frequency forcing regime, we demonstrate the existence of the primary maximal canards 
induced by folded saddle-nodes of type I,
and establish explicit formulas
for the parameter values 
at which the primary maximal canards and their fold curves exist.
Then, we turn to the intermediate- and high-frequency forcing regimes, and show that the forced van der Pol possesses torus canards instead. These torus canards consist of long segments near families of attracting and repelling limit cycles of the fast system, in alternation. We also derive explicit formulas for the parameter values at which the maximal torus canards and their fold curves exist. Primary maximal canards and maximal torus canards correspond geometrically to the situation in which the persistent manifolds near the family of attracting limit cycles coincide to all orders with the persistent manifolds that lie near the family of repelling limit cycles. 
The formulas derived for the folds of maximal canards in all three frequency regimes turn out to be representations of a single formula in the appropriate parameter regimes, and this unification confirms the central numerical observation that the folds of the 
maximal canards created in the low-frequency regime continue directly into the fold curves of the maximal torus canards that exist in the intermediate- and high-frequency regimes. In addition, we study the secondary canards induced by the folded singularities in the low-frequency regime and find that the fold curves of the secondary canards turn around in the intermediate-frequency regime, instead of continuing into the high-frequency regime as the primary maximal canards do. Also, we identify the mechanism responsible for this turning.
Finally, we show that the forced van der Pol equation is a normal form type equation for a class of single-frequency periodically-driven slow/fast systems with two fast variables and one slow variable which possess a nondegenerate fold of limit cycles. The analytic techniques used herein rely on geometric desingularization, invariant manifold theory, Melnikov theory, and normal form methods. The numerical methods used herein were developed in \cite{Desroches2008a,Desroches2010}.

\vspace{8pt}\noindent \textbf{Keywords}\quad folded singularities, canards, torus canards, torus bifurcation, mixed-mode oscillations

\vspace{10pt}\noindent \textbf{AMS subject classifications} 
Primary: 34E17, 34E15, 34D26, 70K70; 
Secondary: 34E13, 34D15, 34C29, 34C45, 37G15, 92C20, 70K43
\end{abstract}
%------------------------------------------------------------

%---------------------------------------------------------------------------------
\section{Introduction}	\label{sec:intro}
%---------------------------------------------------------------------------------

The forced van der Pol equation
is a fundamental model for oscillatory processes 
in physics, electronics, biology, neurology, 
sociology, and economics. 
Possessing strong nonlinear damping effects, 
it is the prototype
of a forced relaxation oscillator,
exhibiting slow and fast time scales,
see \cite{Bold2003, Cartwright1945, Cartwright1950, 
Flaherty1978, Haiduc2009,
Han2012, Levi1981, Levinson1949, 
Sekikawa2005, 
vanderPol1920, vanderPol1927}. 
The equations may be formulated as
\begin{equation}	
\label{fvdp}
\begin{split}
x' &= y-f(x), \\
y' &= \eps( -x + a+b \cos \theta), \\
\theta' &= \omega,
\end{split}
\end{equation}
where the prime denotes the derivative with respect to the fast time variable $\tau$, $f(x)=\frac{1}{3}x^3-x$, and $0 < \eps \ll 1$. The external signal, $a+b \cos \theta$, models a time-periodic driving force with drive frequency $\omega>0$. We refer to equation \eqref{fvdp} as the fvdP equation. Throughout this article, we will work with the form of the system given by \eqref{fvdp}, as it allows us to explore the full range of forcing frequencies $\omega > 0$.

A number of detailed studies of forced van der Pol equations have been carried out in the low-frequency forcing regime, $\omega = \mathcal{O}(\eps)$, see \cite{Bold2003,Guckenheimer2003,Han2012,Szmolyan2004}. We cite in particular the study \cite{Bold2003} of a forced van der Pol system with low-frequency forcing, which presents a detailed analysis of the folded saddle singularities and their attendant canards. 
In the context of excitable systems (in particular, neuronal models), the folded saddle maximal canard plays the role of an excitability threshold manifold, locally dividing trajectories between those that jump at the fold to a different attracting manifold and those that do not \cite{Mitry2013,Wechselberger2013}. This is also true in planar neuronal systems where solutions containing maximal canard segments correspond to excitability thresholds both in the case of type I neurons (integrators) and type II neurons (resonators) \cite{Desroches2013,Izhikevich2000}.
%The folded saddle maximal canard plays the role of a firing threshold manifold, locally dividing trajectories between those that jump at the fold to a different attracting manifold and those that do not. The folded saddle maximal canard plays a similar role in a model of paradoxical excitation in propofol anaesthesia \cite{Mitry2013} and in non-autonomous excitability models \cite{Wechselberger2013}. 
More generally, the canards induced by folded singularities, of folded node, folded saddle, and folded saddle-node types, have also been studied in models of neuronal dynamics \cite{Rotstein2008,Rubin2008,Teka2011}
and in many other systems, see for example \cite{Benoit1990,Brons2006,Desroches2008a,Desroches2012,Krupa2010,
Vo2014,Wechselberger2005,Wechselberger2012,Wechselberger2013}.

In this article, we examine the fvdP equation \eqref{fvdp} in three different regimes of forcing frequencies: low-frequency ($\omega={\cal O}(\eps)$), intermediate-frequency ($\omega={\cal O}(\sqrt{\eps})$), and high-frequency ($\omega={\cal O}(1)$). In each regime, we study the canard solutions that the fvdP equation \eqref{fvdp} exhibits.

%We begin in the low-frequency regime. First, we briefly apply the theory of folded singularities to \eqref{fvdp}, to identify the different types of primary and secondary canards that exist. 
%Then, we place special emphasis on folded saddle-node points of type I (FSN I), and we demonstrate the existence of the primary strong canards. In particular, the following theorem is the first main result of this article:
%In particular, we place special emphasis on folded saddle-node points of type I (FSN I), a close analysis of which leads to the following first main result of this article:
We begin in the low-frequency regime. First, we briefly apply the theory
of folded singularities to \eqref{fvdp}, to identify the different types
of folded singularities that it exhibits in this regime.  We place special
emphasis on folded saddle-node singularities of type I (FSN I), which are
known to generate a number of different types of canard solutions.

The graph of the fast null-cline, $y=f(x)$, of system \eqref{fvdp} with
$\eps=0$ plays a central role in understanding the system dynamics. We are
especially interested in the repelling branch in the middle and the
attracting branch on the right, respectively, of the graph. Let $S_r$
denote the two-dimensional manifold formed by rotating the (middle)
repelling branch through one complete revolution in the angle $\theta:
[0,2\pi)$, and similarly let $S_a$ be the the two-dimensional manifold
formed by rotating the (right) attracting branch through one complete
revolution in $\theta$. In the low-frequency regime of \eqref{fvdp},
Fenichel theory \cite{Fenichel1979,Jones1995} guarantees that, when $\eps$
is sufficiently small, there exist two-dimensional, locally-invariant manifolds
$S_r^{\eps}$ and $S_a^{\eps}$ near $S_r$ and $S_a$, respectively, away
from the fold regions.  In the low-frequency forcing regime of
\eqref{fvdp}, these persistent manifolds are referred to as slow
manifolds, since the dynamics on them is slow in $y$ and $\theta$.

\begin{figure}[ht]
\centering
\includegraphics[scale=0.5]{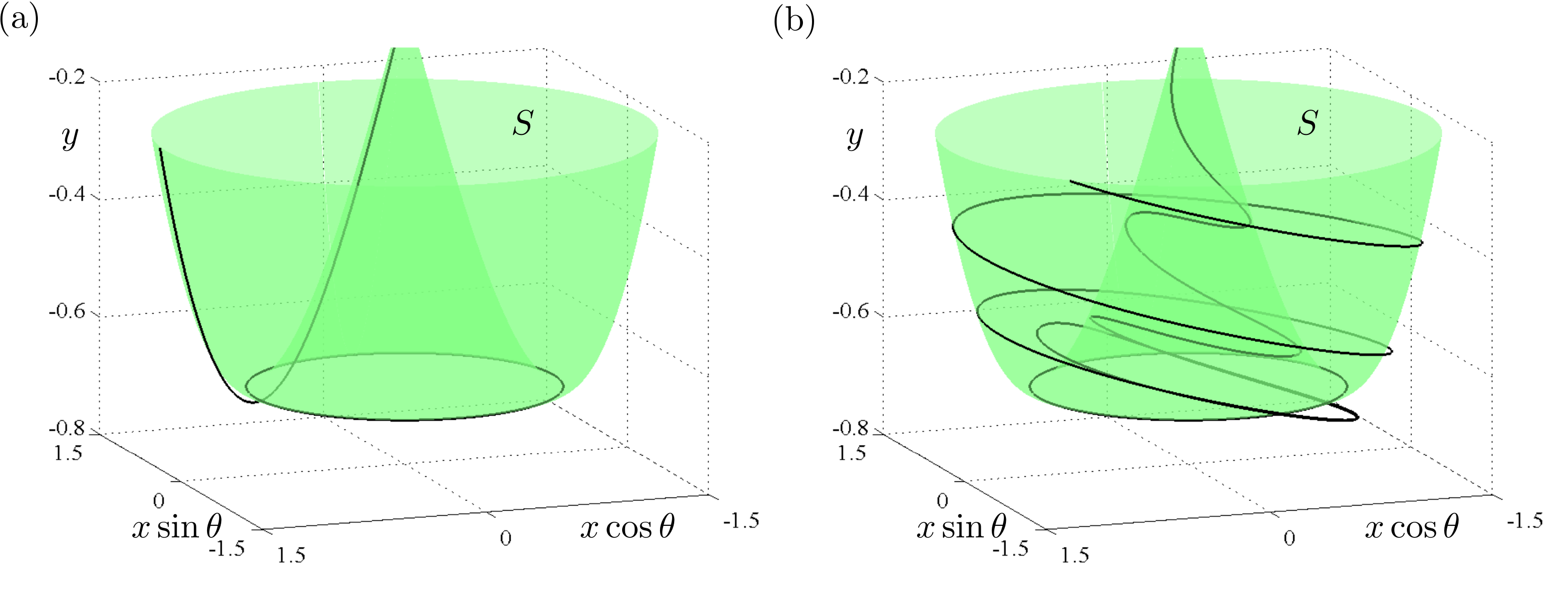}
\vspace*{-0.2cm}
\caption{\small{(a) Segment of a primary maximal canard solution of \eqref{fvdp}, and (b) segment of a maximal torus canard solution of \eqref{fvdp}. Both have long segments near the family of attracting limit cycles (outer portion of the green surface) and near the family of repelling limit cycles (inner portion of the green surface). Here, $a=0.997, b=0.994, \eps=0.02$, and (a) $\omega = 0.001$ and (b) $\omega=0.3$.}}
\label{fig:canb}
\end{figure}

The primary canards of folded singularities are orbits that have a long segment
close to $S_a^{\eps}$, pass through a neighborhood of the folded singularity, and
then have a long segment near $S_r^{\eps}$. The lengths of these segments
depend on the parameter values;  and, there are curves of parameter
values along which the segment near $S_r^{\eps}$ has maximal length, going
all the way up to the other fold curve. These primary canards are referred
to as primary maximal canards. See Figure~\ref{fig:canb}(a) for a
representative primary maximal canard. They are determined geometrically
by the parameter values for which $S_r^{\eps}$ and $S_a^{\eps}$ agree to all
orders in $\eps$, in a manner that is analogous to the maximal limit cycle
canards (also known as the maximal headless ducks) in the classical, planar
van der Pol equation, recall \cite{Benoit1981, Diener1984,Eckhaus1983}.

The following is the first main result of this article:
\begin{theorem}
\label{thm:lowfrequency}
{\bf low-frequency forcing.}
Let $\omega=\eps \overline{\omega}$, where $\overline{\omega}={\cal O}(1)$, and let $b=\mathcal{O}(\sqrt{\eps})$. Then, there exists an $\eps_0>0$ such that for all $0<\eps<\eps_0$, there are two curves in the $(a,\overline{\omega})$ parameter plane given by
\begin{equation}	\label{mainformulas-Reg1}
a = 1 -\frac{\eps}{8} \pm b\exp\left(-\frac{\eps\overline{\omega}^2}{2} \right),
\end{equation}
emanating from the points $(\overline{\omega},a) = (0,1-\frac{\eps}{8} \pm b)$, along which the system \eqref{fvdp} has folds of primary maximal canards. Moreover, for each ${\cal O}(1)$ value of $\overline{\omega}$, the system \eqref{fvdp} has two primary maximal canards for every value of $a$ in the interval between the points on these fold curves. Finally, there are no primary maximal canards for values of $a$ outside the closures of these intervals.
\end{theorem}
This first theorem is established by using the geometric desingularization method, also known as the blow-up method \cite{Dumortier1996,Dumortier2001,Krupa2001}, to inflate the folded saddle node points of type I into hyperspheres, and then by employing invariant manifold theory and Melnikov theory in the appropriate coordinate charts.

Next, we show that the fvdP system \eqref{fvdp} has torus canards in both the regime of intermediate-frequency forcing, in which the system \eqref{fvdp} has three time scales: $x$ is a fast variable, $\theta$ is an intermediate time variable, and $y$ is a slow variable;
and in the regime of high-frequency forcing in which \eqref{fvdp} is a two-fast ($x,\theta$) and one-slow system ($y$). Torus canards are a relatively new type of canard solution discovered in a model of neuronal activity in Purkinje cells \cite{Kramer2008}. 
They consist of long segments near families of attracting and repelling limit cycles of the fast system, in alternation.
Torus canards have recently been shown to exist in a broad array of models, including in three models of neuronal bursting, see \cite{Burke2012}: the Hindmarsh-Rose model, the Morris-Lecar-Terman system, and the Wilson-Cowan-Izhikevich equations; in a model of elliptic bursters,
where the torus canards are rotated versions of limit cycle canards of a planar system, see \cite{Izhikevich2001}; in a rotated van der Pol-type model system, see \cite{Benes2011}; as well as more recently in a model of respiratory rhythm generation in a pre-B\"otzinger complex, see \cite{Roberts2014}. The significance of torus canards is that they play a central role in the transition between periodic spiking and bursting, of different types, in these neuronal models.

%Figure 1: 
%\textcolor{blue}
%{Insert Figure 1 here. This should show 
%a canard for $b={\cal O}(1)$
%together with the critical manifold $S^0$.
%The latex code for including this figure
%is right below these lines, percented out.
%The label for this figure is fig:canb
%Also, note that this was formerly Fig 10 in the 2/19 version.}
%
%HERE
%
%THERE

In the intermediate- and high-frequency regimes of \eqref{fvdp}, Fenichel
theory \cite{Fenichel1979,Jones1995} also guarantees that, when $\eps$ is sufficiently
small, there exist two-dimensional, locally-invariant manifolds near
$S_r$ and $S_a$, away from the fold regions.  We again denote these by
$S_r^{\eps}$ and $S_a^{\eps}$ and label them as persistent manifolds.
However, it is crucial to observe that these persistent manifolds are no
longer slow manifolds in these regions.  Instead, the orbits of
\eqref{fvdp} on these persistent manifolds exhibit two time scales, with
fast rotation due to the limit cycles, as well as slow drift in the
vertical direction, down $S_a^{\eps}$ and up along $S_r^{\eps}$.

Torus canards are orbits of \eqref{fvdp} in the
intermediate- and high-frequency regimes that have long segments near
$S_a^{\eps}$, spiral through a neighborhood of the fold curve of limit
cycles, and then have a long segment near $S_r^{\eps}$. The lengths of time
that torus canards spiral around near $S_a^{\eps}$ and $S_r^{\eps}$ depend on
the system parameters, and for system \eqref{fvdp} there are curves of
parameter values along which the time spent near
$S_r^{\eps}$ is maximal, with the orbits spiralling all the way up
$S_r^{\eps}$. These are defined to be maximal torus canards, in analogy with
the maximal limit cycle canards of the unforced van der Pol oscillator. A
representative maximal torus canard is shown in Figure~\ref{fig:canb}(b).

For system \eqref{fvdp} in the intermediate-frequency regime, we prove the following theorem:
\begin{theorem}
\label{thm:intermediatefrequency}
{\bf intermediate-frequency forcing.}
Let $\omega=\sqrt{\eps} \Omega$, where $\Omega={\cal O}(1)$, and let $b=\mathcal{O}(\eps)$. Then, there exists an $\eps_0>0$ such that for all $0<\eps<\eps_0$, there are two curves in the $(a,\Omega)$ parameter plane given by
\begin{equation}	\label{mainformulas-Reg2}
a = 1 - \frac{\eps}{8} \pm b \exp \left( - \frac{\Omega^2}{2} \right).
\end{equation}
along which the system \eqref{fvdp} has folds of maximal canards. Moreover, for each fixed ${\cal O}(1)$ value of $\Omega$, the system \eqref{fvdp} has two canards for every value of $a$ in the interval between these fold curves, and none outside the closure of these intervals.
\end{theorem}
%
%The folds of torus canards may also be referred to as maximal torus canards. 
This theorem is also established by using the geometric desingularization method; however, in this regime, we inflate the circular fold curve along which the attracting and repelling limit cycles meet into a two-torus, rather than the FSN I points. Also, the scalings are different, as is the analysis in the coordinate charts near the torus.

Then, for the high-frequency regime, we establish:
\begin{corollary}
\label{cor:highfrequency}
{\bf high-frequency forcing.}
Let $\omega={\cal O}(1)$, and let $b=\mathcal{O}(\eps)$. Then, there exists an $\eps_0>0$ such that for all $0<\eps<\eps_0$, there are two curves in the $(a,\omega)$ parameter plane given by
\begin{equation}	\label{mainformulas-Reg3}
a = 1 - \frac{\eps}{8} \pm b \exp \left( - \frac{\omega^2}{2\eps} \right).
\end{equation}
along which the system \eqref{fvdp} has folds of torus canards. Moreover, there exists a pair of torus canards for each parameter value
in the interval between these fold curves.
\end{corollary}

%The torus canards whose existence is established in Theorem~\ref{thm:intermediatefrequency} and Corollary~\ref{cor:highfrequency}
%are solutions of \eqref{fvdp} that consist of alternating long segments near branches of attracting limit cycles and branches of repelling limit cycles of the fast system. A representative torus canard is shown in Figure \ref{fig:canb}. Also, w
We note that the presence of the torus canards in this type of fast-slow system is signalled by the existence of a fold of limit cycles of the fast system, here at $(x,y)=(1,-\frac{2}{3})$, together with a nearby torus bifurcation in the full system, here at
\begin{equation} \label{eq-TBif}
1 - a^2 - \frac{1}{2} \frac{b^2 \eps^2}{ (a^2-1)^2\omega^2 + (\eps-\omega^2)^2} = 0.
\end{equation}
See Appendix A. These two triggering mechanisms arise ubiquitously in fast-slow systems with two or more fast variables.

Having established these theorems for the existence of the primary maximal canards and the torus canards, as well as their folds,
we now analyze the relationship between these results. Plainly, the formulas for the curves of folds \eqref{mainformulas-Reg1},
\eqref{mainformulas-Reg2}, and \eqref{mainformulas-Reg3} in the three different regimes are all representations of the same formula,
\begin{equation}	\label{mainformulas}
a = 1 - \frac{\eps}{8} \pm b \exp \left( - \frac{\omega^2}{2\eps} \right),
\end{equation}
in the respective frequency regimes.
The exponential term has magnitude $b$ (which is ${\cal O}(\sqrt{\eps}))$ and is slowly varying with $\overline{\omega}$ in the low-frequency regime (Theorem~\ref{thm:lowfrequency}), small-amplitude ($b=\mathcal{O}(\eps)$) and varying with ${\cal O}(1)$ frequency $\Omega$ in the intermediate-frequency regime (Theorem~\ref{thm:intermediatefrequency}), and exponentially small in $\eps$ in the high-frequency regime (Corollary~\ref{cor:highfrequency}).

%FIG 2 here
%for the two curves of folds of canards in the $\omega-a$ plane
%for the same three representative values of $b$.
%The label for this figure is fig:agreement
%At the moment there is only one frame here, ie one value of $b$}
%HERE
%
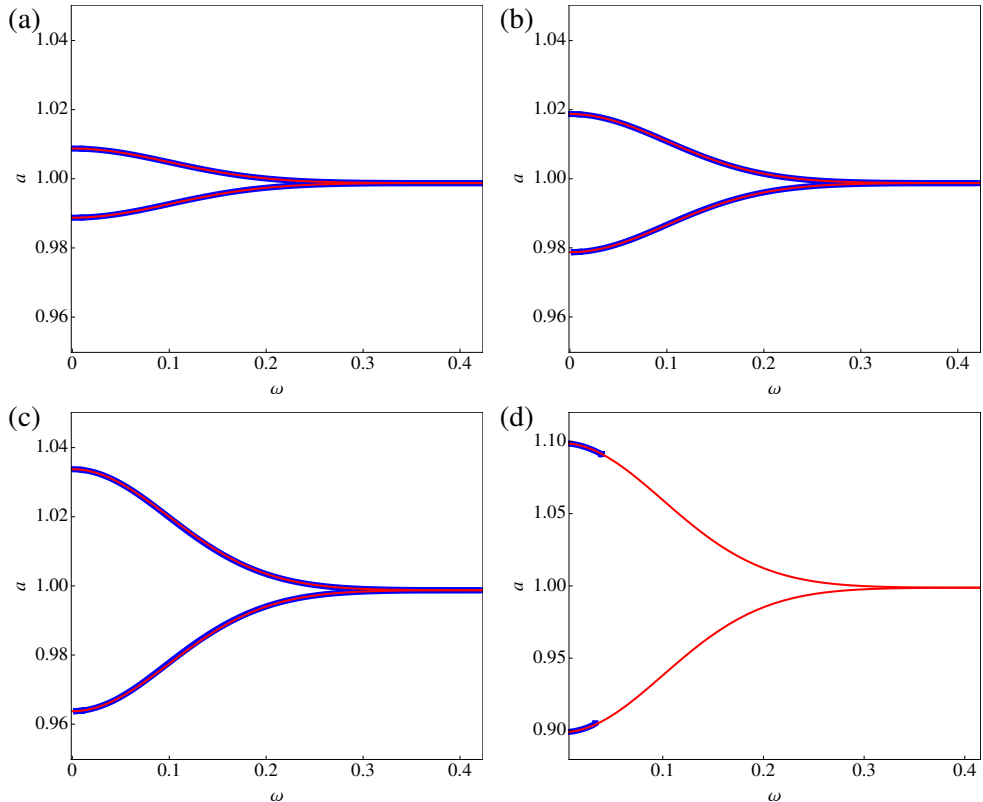
\begin{figure}[!t]
\centering
\includegraphics[width=5in]{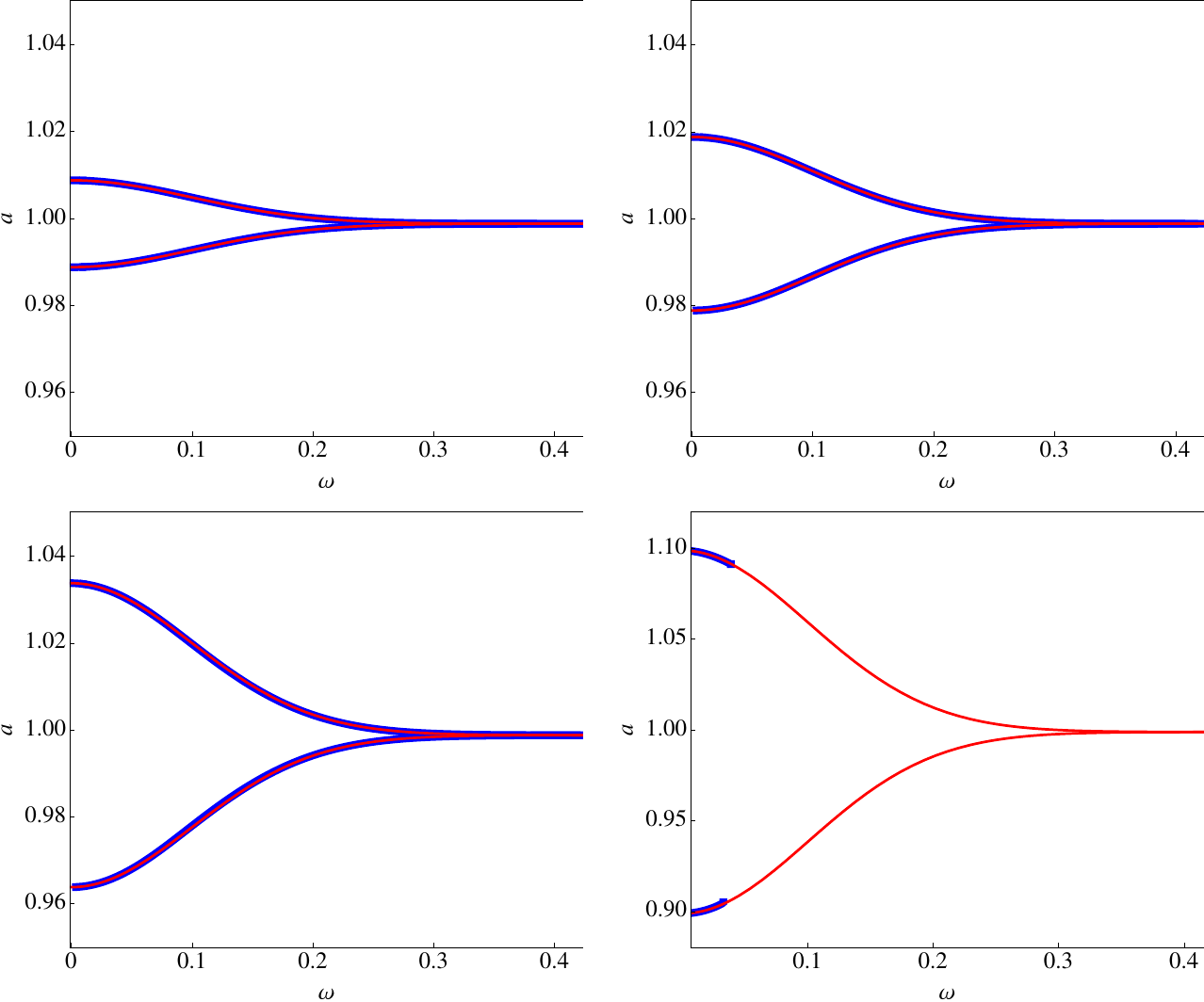} 
\put(-364,292){(a)}
\put(-180,292){(b)}
\put(-364,142){(c)}
\put(-180,142){(d)}
\vspace*{-0.2cm}
\caption{Curves of folds of maximal canards in the $(\omega,a)$ plane as obtained from \eqref{mainformulas} (red curves) and 2-parameter numerical continuation (blue curves) for $\eps=0.01$ and (a) $b=0.01$; (b) $b= 0.02$; (c) $b= 0.035$l and, (d) $b=0.1$. For $b = \mathcal{O}(\eps)$ (panels (a)--(c)), there is good agreement between theoretical (red) and numerical (blue) results over the entire range of forcing frequencies, including for both the primary maximal canards which exist for $\omega={\cal O}(\eps)$ and the maximal torus canards which exist for $\omega={\cal O}(1)$. Note that the scales in panels (a)--(c) are the same. For $b=\mathcal{O}(\sqrt{\eps})$ (panel (d), in which the vertical scale is different), we find that the numerical continuation terminates when $\omega$ is no longer  $\mathcal{O}(\eps)$.  
}
\label{fig:agreement}
\end{figure}
%
%THERE

The analysis in all three regions shows that the values of the parameter $a$ for which the canards exist in between the fold curves may similarly be summarized succinctly in one formula:
\begin{equation}	\label{theta0-mainformulas}
a = 1 - \frac{\eps}{8}- b\cos(\theta_0)\exp \left(-\frac{\omega^2}{2\eps}\right).
\end{equation}
Here, $\theta_0$ is an arbitrary phase, and the magnitude and dynamics of the exponential term are also as discussed above.

It is also of interest to observe that,
in the limit of 
$\omega \to \infty$,
formulas \eqref{mainformulas} and \eqref{theta0-mainformulas}
become
$a \to a_c := 1-\frac{\eps}{8}$, 
which corresponds exactly 
to the leading order locations 
of the maximal limit cycle canards
in the planar vdP equation,
see for example \cite{Krupa2001}. 
Therefore, as expected,
for sufficiently high-frequency forcing, 
the effect of the forcing averages out 
to this order, and \eqref{fvdp} 
behaves like the classical planar vdP equation.
In this limiting regime,
the torus canards of \eqref{fvdp} appear to be rotated copies of the limit cycle canards of the planar vdP.

Then, with the above analytical results in hand,
we turn next to the results of numerical continuations
which confirm that 
the curves of the folds of primary strong canards 
observed in the low-frequency regime
continue directly to the curves of the folds of torus canards
discovered in the intermediate- and high-frequency regimes. 
This is illustrated in Figure~\ref{fig:agreement}.
Moreover, as is also shown in this figure,
the agreement between the formulas and the numerical continuation
results are excellent within the parameter regions
stated in the theorems. 
We also note that the theory does not appear to extend
outside of these regions,
and also preliminary results of numerical continuation
reveal different dynamics there.
Overall, then, the formulas 
\eqref{mainformulas},
and \eqref{theta0-mainformulas}
together with the numerical continuations,
will directly imply that
the primary strong canards,
which exist in the low-frequency forcing region,
continue naturally to the branches of torus canards,
which exist in the high-frequency regime,
where the folded singularities cease to exist,
with the transition happening 
in the intermediate-frequency regime.
Understanding the continuation dynamics
of these curves is one of the main results of this article.
Moreover, the results here will also help shed light
on other models with torus canards.
In particular, we observe that numerical simulations 
of a rotated van der Pol-type model 
exhibit the same continuation of the maximal canards
across the entire range of forcing frequencies;
see Figure 5 in \cite{Benes2011}.
Numerical continuations in other neuronal (or neural) models \cite{Burke2012}
show similar phenomena.

\begin{figure}[!t]
\centering
\includegraphics[scale=0.5]{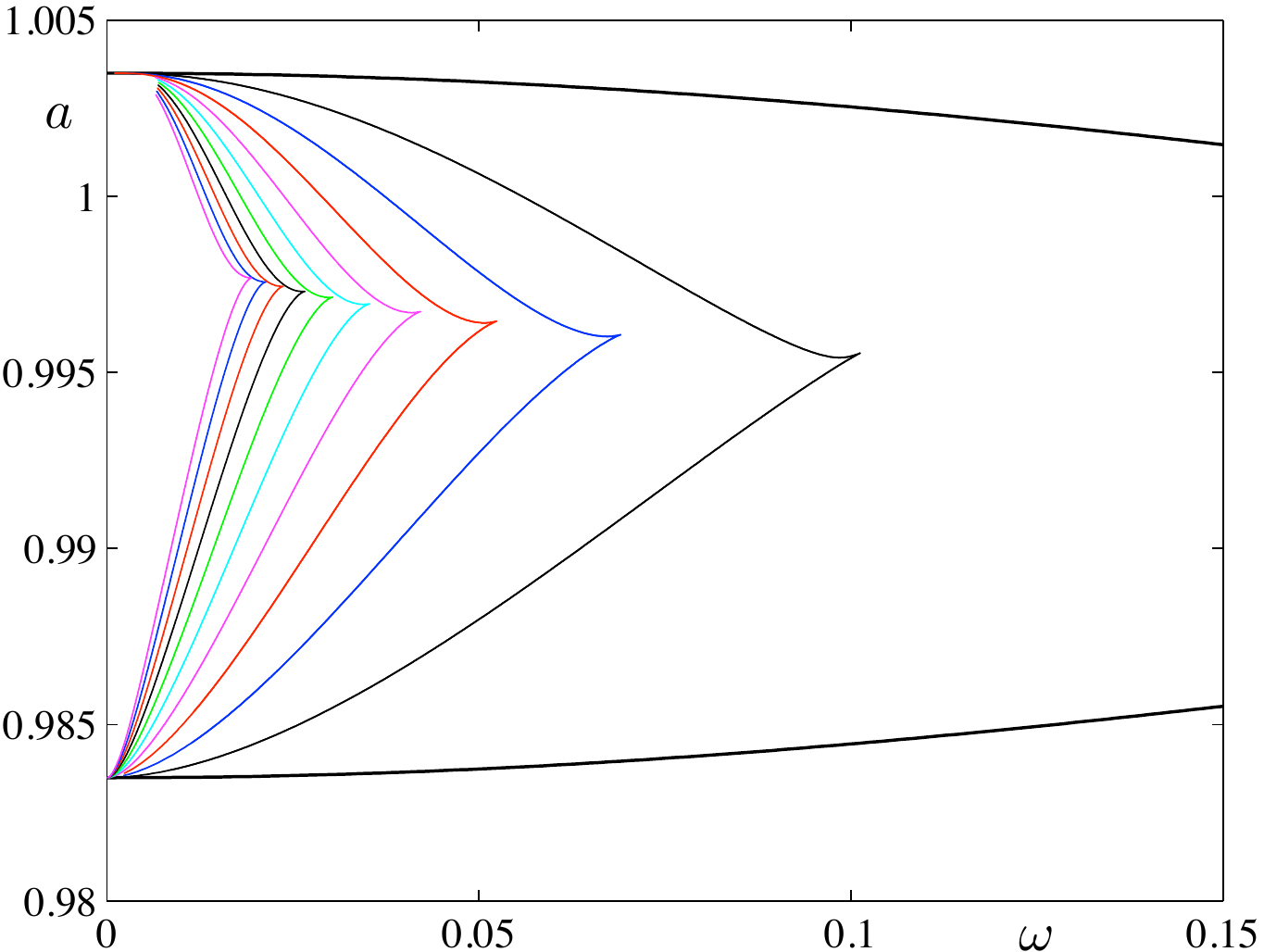}
\vspace*{-0.2cm}
\caption{Curves of folds along the computed branches of primary and secondary canards. Here, $\eps=0.05$ and $b=0.01$.}
\label{fig:foldcancont}
\end{figure}

In this article, we also study the secondary canards of \eqref{fvdp}.
Secondary canards, lie near the primary canards for most of their lengths,
and make a number of small loops around an axis of rotation usually referred to as the \textit{weak canard}.
Secondary canards are indexed by the number of loops they make around the weak canard
and by the value of the $y$-intercepts
of the solutions during
the nearly-horizontal jumps
that occur from a neighborhood of the family of repelling limit cycles
back to a neighborhood of the family of attracting limit cycles.
In particular, we carry out numerical continuation
of the branches of maximal secondary canards of \eqref{fvdp}
that are created by the folded singularities
in the low-frequency regime.
In contrast to what we find for the primary canards,
the branches of the secondary canards
turn around well before they reach the high-frequency regime.
See Figure~\ref{fig:foldcancont}.
Also, we identify the mechanisms
which cause the branches to turn.

%Figure 3: 
%\textcolor{blue}
%{Insert Figure 3 here. This should show 
%the curves of the folds of secondary canards in the $\omega-a$ plane
%for one value of $b$.
%The latex code for including this figure
%is right below these lines, percented out.
%The label for this figure is fig:foldcancont
%This was formerly Fig 11 in the 2/19 version.}
%HERE
%
%
%THERE

To conclude this article, we demonstrate that \eqref{fvdp} serves as a local normal form for slow/fast systems with one slow variable
and two fast variables in which the fast subsystem possesses a non-degenerate fold of limit cycles, and in which the slow system is subject to time-periodic forcing. These fast-slow systems exhibit torus canard explosions, just as shown here for the fvdP equation \eqref{fvdp}, and one may therefore directly conclude, by applying the same techniques used herein, that the folds of their canards
behave in a similar fashion.

%Figure 4.
%\textcolor{blue}{ formerly figure 2 in the 2/19 version}
%HERE
%
\begin{figure}[h]
\centering
\includegraphics[width=5in]{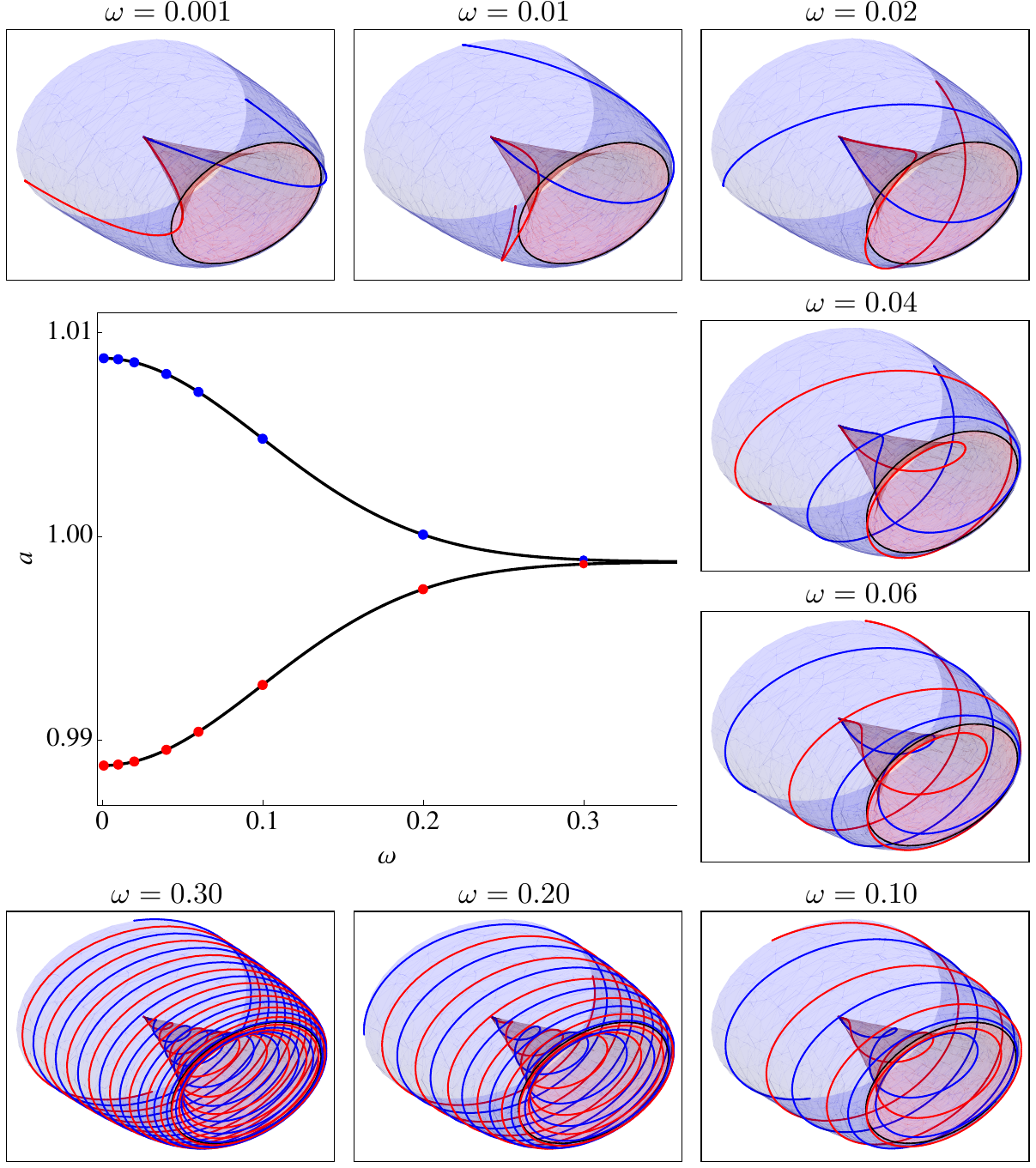}
\caption{2-parameter continuation of folds of primary maximal canards of a folded saddle-node (type I) for $\eps = 0.01$ and $b=0.01$.
Orbit segments are plotted in `Cartesian' coordinates $(u,v,y) = (x \cos \theta, x \sin \theta, y)$. }
\label{fig:motivation}
\end{figure}
%
%THERE

The numerical method 
developed in \cite{Desroches2008a, Desroches2010},
and also employed in \cite{Desroches2012},
is the main numerical method used throughout this article
to find the persistent
invariant manifolds and the curves of maximal canards
that lie in the intersections of these manifolds.
This method,
which uses the AUTO continuation software
\cite{Doedel2007},
turns the problem of finding the invariant
manifolds of fast-slow systems into a boundary value problem
for system \eqref{fvdp} with the time of integration
included as a parameter.
Then, the parametrized families of solutions of the two-point boundary
value problems are continued.
This allows to integrate in positive and negative time 
using pseudo-arclength continuation, 
approximating the orbit segments of solutions
of system~\eqref{fvdp} subject to particular boundary conditions 
by orthogonal collocation, 
which is very well suited to multiple time scale vector fields (see \cite{Desroches2012,Vo2014}).

The outline of the article is as follows. 
In Section \ref{sec:foldedsingularities}, 
we consider system \eqref{fvdp} with low-frequency forcing,
$\omega = \mathcal{O}(\eps)$,
and apply the theory of folded singularities
in a straightforward manner 
to find the associated primary and secondary canards
of folded singularities.
Then, in Section~\ref{sec:AnalysisRegion1}, 
we prove Theorem~\ref{thm:lowfrequency},
establishing the existence of the primary maximal canards
induced by the FSN I points,
and the associated fold curves,
including the derivations of formulas 
\eqref{mainformulas-Reg1} 
and \eqref{theta0-mainformulas} 
for system \eqref{fvdp} in 
the low-frequency 
${\cal O}(\eps)$ region.
We then turn to the cases of
intermediate- and
high-frequency forcing
in Section \ref{sec:AnalysisRegion2+3}, 
where we study the torus canards of \eqref{fvdp}.
We prove Theorem~\ref{thm:intermediatefrequency}
and Corollary~\ref{cor:highfrequency},
establishing the existence of the maximal torus canards and their fold curves
in these regimes, as given by the formulas 
\eqref{mainformulas-Reg2}
and \eqref{mainformulas-Reg3}.
This shows analytically that the curves of the fold curves of the primary maximal 
canards, which are born in the low-frequency regime,
continue for all $\omega>0$ into the fold curves of the maximal torus canards.
We also observe that the analytically-derived formulas 
and the curves obtained
in the numerical continuations
agree over the entire range of forcing frequencies.
Then, in Section~\ref{sec:seccan}, 
we examine the numerical continuations
of the folds of secondary canards,
and we identify the mechanism by which they turn around
well before they reach the high-frequency regime.
Also, we analyse how the curves
of the folds of secondary canards 
induced by folded nodes change
as the parameter $b$ is varied, up to and including $b=\mathcal{O}(1)$. 
and hence as the distance between
the folded node and the folded saddle
is varied.
The final main result
of this article is presented in Section~\ref{sec:normalform}.
We demonstrate that \eqref{fvdp} 
may be considered as a local normal form 
for some generic fast-slow systems that have a fold of limit cycles
and that undergo a torus canard explosion.
In Appendix A, we prove, using second-order averaging, the existence of a torus bifurcation in \eqref{fvdp} and calculate the locus \eqref{eq-TBif} for this bifurcation in parameter space.

%---------------------------------------------------------------------------------
\section{Low-Frequency Forcing: Canards of Folded Singularities, Especially of FSN I Points}	\label{sec:foldedsingularities}	
%---------------------------------------------------------------------------------

In this section, we present a brief review and analysis 
of the folded singularities that system \eqref{fvdp} possesses 
in the regime of low drive frequency, 
i.e. $\omega = \eps \overline{\omega}$, 
where $\overline{\omega} = \mathcal{O}(1)$. 
Readers familiar with the theory of folded singularities
and their canards may proceed to Section~\ref{sec:AnalysisRegion1}.
In this regime, \eqref{fvdp} is
\begin{equation}	\label{fvdpfasttime}
\begin{split}
x^\prime &= y-f(x),\\ 
y^\prime &= \eps \left( -x + a+b \cos \theta \right), \\
\theta^\prime &= \eps \overline{\omega}.
\end{split}
\end{equation}
It is a 1-fast/2-slow problem with fast variable $x$
and slow variables $(y,\theta)$. 
We analyse the reduced dynamics 
associated to \eqref{fvdpfasttime} 
and derive the desingularised vector field 
on the critical manifold. 
Then, we identify the canards of the folded singularities.

%------------------------------------------
\subsection{The Layer and Reduced Systems} 
%------------------------------------------

Taking the singular limit $\eps \to 0$ in \eqref{fvdpfasttime}, 
one finds the 1D layer problem
\begin{equation}	\label{fvdplayer}
\begin{split}
x^\prime &= y-f(x),
\end{split}
\end{equation}
where $y$ and $\theta$ are parameters. 
Alternatively, the singular limit $\eps \to 0$ in \eqref{fvdp} 
gives the 2D reduced system
\begin{equation}	\label{fvdpreduced}
\begin{split}
0 &= y-f(x),\\ 
\dot{y} &= -x + a+b \cos \theta, \\
\dot{\theta} &= \overline{\omega},
\end{split}
\end{equation}
where the overdot denotes the derivative with respect to the slow time 
$t= \eps \tau$.
The manifolds $S_a^{\eps}$ and $S_r^{\eps}$ are non-unique.
Hence, the canards that lie near the manifolds
are also non-unique.
However, for a fixed choice of invariant manifolds, $S_a^{\varepsilon}$ and $S_r^{\varepsilon}$, their transverse intersections correspond to maximal canards.

Systems \eqref{fvdplayer} and \eqref{fvdpreduced} 
are two different approximations of the forced vdP equation. 
The idea of 
Geometric Singular Perturbation Theory (GSPT) \cite{Jones1995} 
is to combine information from the 1D layer 
and 2D reduced problems in order to understand the dynamics 
of the full 3D forced vdP equation for $0<\eps \ll 1$.

We begin with an analysis of the 1D layer problem \eqref{fvdplayer}, 
which is an approximation of \eqref{fvdp} 
wherein the slow processes are assumed to move so slowly 
that they are fixed. 
The \emph{critical manifold} 
is the set of equilibria of the layer problem \eqref{fvdplayer}:
\[ S := \left\{ (x,y,\theta) \in \mathbb{R}^2 \times \mathbb{S}^1 : y = f(x) \right\}. \]
Linear stability analysis of \eqref{fvdplayer} 
shows that there are disjoint curves, $L$, of fold points given by
\[ L := \left\{ (x,y,\theta) \in S : x = \pm 1 \right\}, \]
which separate the outer attracting sheets, $S_a$, of $S$ from the middle repelling sheet, $S_r$, of $S$. Fenichel theory \cite{Fenichel1979,Jones1995} guarantees that the normally hyperbolic segments of $S$ (i.e. the parts of $S_a$ and $S_r$ at $\mathcal{O}(1)$ distances from the fold curve $L$) will persist as invariant slow manifolds, $S_a^{\eps}$ and $S_r^{\eps}$, of \eqref{fvdp} for $0<\eps \ll 1$. 

The price we pay for the approximation \eqref{fvdplayer} is that we have trivial dynamics on $S$. To obtain a non-trivial flow on the critical manifold, we turn to the reduced problem. The 2D reduced problem \eqref{fvdpreduced} is an approximation of \eqref{fvdp} wherein the fast motions are assumed to be so rapid that they immediately settle to their steady state, which is precisely the critical manifold. In other words, the reduced problem prescribes a non-trivial flow along $S$. The price we pay for this particular approximation is that the reduced flow is not defined away from the critical manifold. Note that the restriction of the flow of \eqref{fvdp} to $S^{\eps}$ is a small smooth perturbation of the reduced flow on $S$.

To analyse the flow on a manifold, we use the coordinates
$(x,\theta)$. 
We differentiate the algebraic constraint $y=f(x)$ with respect to $t$
to obtain the evolution equations in this coordinate chart,
\begin{equation}	\label{reducedprojection}
\begin{split}
(x^2-1) \dot{x} &= -x+a+b \cos \theta, \\
\dot{\theta} &= \overline{\omega}.
\end{split}
\end{equation}
The reduced flow \eqref{reducedprojection} is singular along the fold points $L$ of \eqref{fvdplayer}. 
To remove the finite time blow-up of solutions at the folds, we rescale time $dt = (x^2-1)\, ds$ to obtain 
\begin{equation}	\label{desing}
\begin{split}
\dot{x} &= -x+a+b \cos \theta, \\
\dot{\theta} &= \overline{\omega} (x^2-1),
\end{split}
\end{equation}
where the overdot now denotes derivatives with respect to $s$. 
System \eqref{desing} is equivalent to the reduced flow 
\eqref{reducedprojection} 
on the attracting sheets $S_a$, 
where the time rescaling $dt=(x^2-1)\,ds$ 
preserves the orientation of trajectories. 
On the repelling sheet $S_r$, however, we have $x^2-1<0$,
so that the time rescaling reverses the orientation of trajectories. 
Thus, to obtain the reduced flow \eqref{reducedprojection} 
on $S_r$ from \eqref{desing}, 
we simply reverse the direction of trajectories 
of \eqref{desing} 
whenever we are on the repelling sheet of the critical manifold.

%------------------------------------------
\subsection{Folded Singularities \& Singular Canards} 
%------------------------------------------

The desingularised system \eqref{desing} possesses special equilibria called \emph{folded singularities}, $M$, which are points along the fold curves where the right hand side of the $x$-equation vanishes. In system \eqref{desing}, there are infinitely many pairs of such points (when $\theta$ is considered in its lift to $\mathbb{R}$):
\[ M := \left\{ (x,y,\theta) \in L : \theta = 2k \pi \pm \cos^{-1} \left( \frac{1-a}{b} \right), \, k \in \mathbb{Z} \right\}, \]
where $|1-a| \leq b$. Folded singularities are not true equilibria of the reduced flow \eqref{reducedprojection}. Instead, they correspond to points of \eqref{reducedprojection} where there is potentially a cancellation of a simple zero in the $x$-equation and trajectories may pass through the fold (via the folded singularity) with finite speed. Such a trajectory of the reduced flow that passes through a folded singularity and crosses from the attracting (resp. repelling) sheet to the repelling (resp. attracting) sheet is called a \textit{singular canard} (resp. \textit{singular faux canard}) \cite{Szmolyan2001,Wechselberger2005,Wechselberger2012}. 

Considered as equilibria of the desingularised system \eqref{desing}, folded singularities are classified according to their linearization. Folded nodes have real eigenvalues of the same sign. Folded saddles have real eigenvalues of opposite sign, whilst folded foci have complex eigenvalues. In the forced vdP equation \eqref{fvdp}, 
we find that for $\overline{\omega} >0$ the folded singularities with 
\[ \theta_s(k) = 2k\pi - \cos^{-1}\left( \frac{1-a}{b} \right), \] 
are folded saddles, whilst the folded singularities with
\[ \theta_n(k) = 2k\pi + \cos^{-1}\left( \frac{1-a}{b} \right),\] 
are folded nodes provided
\[ (1-a)^2 < b^2 < (1-a)^2+\frac{1}{64\,\overline{\omega}^2}. \]
Folded nodes and folded saddles have been demonstrated to be the organising centres for complex phenomena. Folded nodes for instance have been identified as the cause of the small oscillations in MMOs in various neurophysiological problems \cite{Erchova2008}, such as in a self-coupled FitzHugh-Nagumo model, in a Hodgkin-Huxley model \cite{Rubin2008} and in a pituitary lactotroph cell model \cite{Teka2011}. More recently, folded saddles have been identified as playing a significant role in distinguishing between transient spiking and quiescence in a model of propofol anaesthesia \cite{Mitry2013} and in non-autonomous excitability models \cite{Wechselberger2013}.

%------------------------------------------
\subsection{Canards of Folded Saddle-Node Type I Points} 
\label{sec:2.3}
%------------------------------------------ 

Folded nodes and folded saddles 
can be created through bifurcations in two distinct ways in \eqref{fvdp}: 
via a folded saddle-node (FSN) of type I \cite{Vo2014,Krupa2010} 
or via an FSN of type II \cite{Krupa2010}. 
Both FSN's correspond to a zero eigenvalue of the folded node 
(or folded saddle). 
The two FSN scenarios are distinguished by their geometry. 
In the FSN I limit, 
the center manifold of the FSN I (in system \eqref{desing}) 
is tangential to the fold curve. 
In the FSN II limit, 
the center manifold of the FSN II point is transverse to the fold curve. 
We focus here on FSN I and refer to the remark below
for FSN II points.

\begin{figure}[h]
\centering
\includegraphics[width=5in]{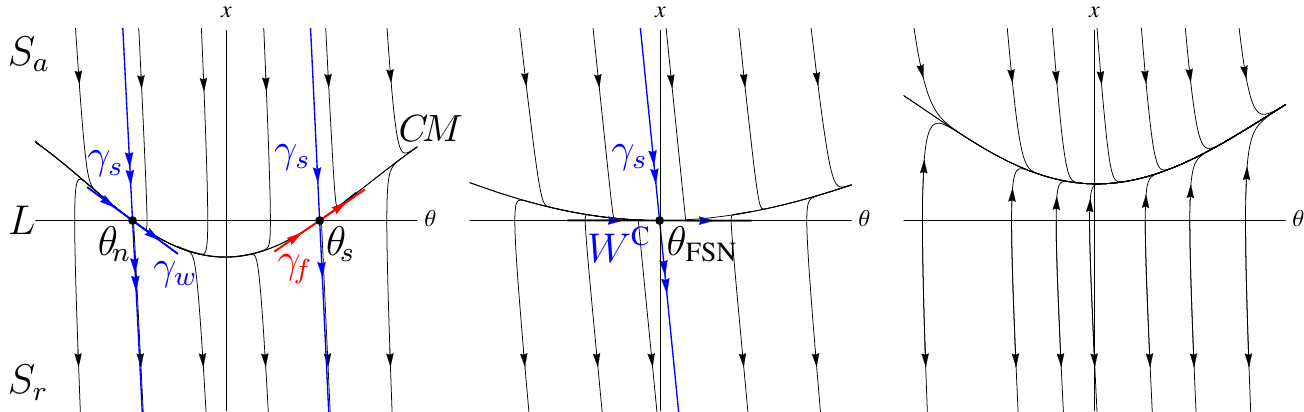}
\put(-364,107){(a)}
\put(-240,107){(b)}
\put(-120,107){(c)}
\caption{\small{Reduced flow \eqref{reducedprojection} of the fvdP equation \eqref{fvdp} shown in a neighbourhood of the upper fold curve $L$ (defined by $x=1,y=-\frac{2}{3}$) for $\overline{\omega}=1, b=0.01$, and (a) $a=1+\frac{b}{2}$, (b) $a=1+b$ and (c) $a=1+\frac{3b}{2}$. In panel (a) where $a<1+b$, there is a folded node ($\theta_n$) and a folded saddle ($\theta_s$). The strong and weak eigendirections of the folded node are denoted by $\gamma_s$ and $\gamma_w$, respectively. Similarly, the singular canard and faux canard of the folded saddle are labelled $\gamma_s$ and $\gamma_f$, respectively. Note that there is a heteroclinic connection $CM$ from $\theta_n$ to $\theta_s$. In particular, $\gamma_w$ is tangent to $CM$ at $\theta_n$, and $\gamma_f$ is tangent to $CM$ at $\theta_s$. In panel (b) where $a=1+b$, the folded singularities merge to a FSN I (indicated by $\theta_{\text{FSN}}$). In this case, the singular strong canard of the folded node merges with the maximal canard of the folded saddle. Meanwhile, the singular weak canard of the folded node merges with the faux canard of the folded saddle to become the center manifold $W^C$ of the FSN I. In panel (c) where $a>1+b$, the folded singularities (and associated singular canards) have been destroyed in the FSN I bifurcation, and there are no canard dynamics.}}
\label{fig:desing}
\end{figure}

The FSN I is the codimension-1 bifurcation 
of the desingularised system \eqref{desing}
in which a folded saddle and a folded node coalesce 
and annihilate each other 
in a saddle-node bifurcation of folded singularities. 
For the forced vdP equation \eqref{fvdp}, 
there are infinitely many such FSN I points:
$(x,y,\theta) = (1,-\frac{2}{3},2k \pi)$,
and they occur for $a=1 \pm b$ and $\omega=\eps\overline{\omega}$.
The FSN I at $a=1-b$ has its center manifold on $S_r$ so that the funnel region (enclosed by the strong canard of the folded node / folded saddle canard and the fold curve) vanishes in the FSN I limit. In this case, we expect generic solutions of \eqref{fvdp} near this FSN I limit to either be relaxation oscillations or MMOs. The FSN I at $a=1+b$ on the other hand has its center manifold on $S_a^+$ so that the funnel persists in the FSN I limit and most solutions of \eqref{fvdp} can tunnel through $S_r$ and return to $S_a$. A representative example of the passage through a FSN I bifurcation at $a=1+b$ is shown in Figure \ref{fig:desing}.

In the FSN limit, the standard folded node/folded saddle theory requires modification. For the FSN I, the following results have recently been proved in \cite{Vo2014}, valid for $0<\eps \ll 1$ and $\mu = \mathcal{O}(\eps^{\alpha})$ where $\alpha \geq 1/4$:
\begin{enumerate}
\item The singular strong canard of the folded node perturbs to the primary maximal strong canard. The singular canard of the folded saddle perturbs to a maximal canard.
\item There exists a heteroclinic connection $CM$ between the folded nodes and folded saddles of \eqref{desing}. This heteroclinic perturbs to a canard-faux canard solution $CM^{\eps}$ that corresponds to both the primary weak canard of the folded node and the faux canard of the folded saddle (faux canards are the equivalent of singular faux canards for $0<\varepsilon\ll 1$). 
\item There exist $\mathcal{O}(\eps^{-1/4})$ canards and faux canards.
\end{enumerate}
Thus, canards and faux canards of the FSN I oscillate about an axis of rotation, which is approximately given by the heteroclinic $CM$ that connects the folded node and folded saddle (see Figure \ref{fig:desing}(a) for instance). 
For the fvdP, we find that 
$CM := \left\{ (x,y,\theta) \in S: x = a+b \cos \theta \right\}. $
We study the associated maximal canards of \eqref{fvdp} in 
Section~\ref{sec:AnalysisRegion1}.

\begin{remark}
For the forced vdP equation \eqref{fvdp}, 
FSN II points are codimension-2 bifurcation points
of the desingularised flow (corresponding to $\overline{\omega} = 0$), 
and constitute a special case of the FSN I. 
They can be analyzed using the approach presented in \cite{Krupa2010}.
\end{remark}

%---------------------------------------------------------------------------------
\section{Loci of the Maximal Canards for Low Forcing Frequencies} 	\label{sec:AnalysisRegion1}
%---------------------------------------------------------------------------------

In this section, we analyse system \eqref{fvdp} with low-frequency forcing ($\omega=\eps\overline{\omega}$). We prove Theorem~\ref{thm:lowfrequency}, demonstrating that, for $b=\mathcal{O}(\sqrt{\eps})$, formula \eqref{mainformulas-Reg1} gives the branches of the folds of the primary maximal canards and that for each value of $a$ in between the fold curves, there are two primary maximal canards of system \eqref{fvdp} given in parameter space by \eqref{theta0-mainformulas}. 
More precisely, we analyze the FSN I points to show that formula \eqref{mainformulas-Reg1} gives the locus of points at which the primary strong canard of a folded node point merges with the folded saddle maximal canard. We present the analysis for the FSN I that occurs for $a=1-b$ and note that the FSN I at $a=1+b$ is treated similarly. 

%set
%\[  a - 1 + b\cos(\theta_0) = \eps\gamma, \]
%where $\gamma$ is ${\cal O}(1)$ with respect to $\eps$. Also, we translate the FSN I to the origin 
For the analysis with low-frequency forcing, we first translate the FSN I at $a=1-b$ to the origin
\[ u= x-1, \quad v=y+\frac{2}{3}, \quad \eta = a-1+b, \]
%
%recall Section~\ref{sec:2.3}.
%Therefore, here,
%the fvdP \eqref{fvdp} is
%conveniently analysed in the form:
%
%\begin{align}   \label{fvdpRegion1}
%u^\prime &= v-u^2-\frac{1}{3} u^3, \nonumber \\
%v^\prime &= \eps \left( - u + \eps\gamma + \sqrt{\eps}\beta \cos(\theta_0) (\cos(\theta) - 1) - \sqrt{\eps}\beta \sin(\theta_0)\sin(\theta)  \right), 
%\nonumber \\
%\theta^\prime &= \eps \overline{\omega}.
%\end{align}
%
%Here, $\theta_0$ is an arbitrary phase.
so that \eqref{fvdp} is transformed to
\begin{equation} \label{fvdpRegion1}
\begin{split}
u^\prime &= v-\left( u^2+\frac{1}{3}u^3 \right), \\
v^\prime &= \eps \left( -u+\eta +b(\cos \theta -1) \right), \\
\theta^\prime &= \eps \overline{\omega}.
\end{split}
\end{equation}
We then inflate the FSN I singularity to a hypersphere using the spherical blow-up transformation: 
\[ u = \overline{r}^2 \overline{u}, \quad v = \overline{r}^4 \overline{v}, \quad \theta = \overline{r} \overline{\theta}, \quad \eps = \overline{r}^4 \overline{\eps}. \]
Moreover, we rescale the parameters $b$ and $\eta$ as
\[ b = \sqrt{\eps} \beta, \quad \eta = \eps \gamma, \]
where $\beta = \mathcal{O}(1)$ and $\gamma = \mathcal{O}(1)$.
Also, we append the trivial equation $\eps'=0$ to system \eqref{fvdpRegion1} and take $\mu>0$ sufficiently small. The (spherical) blow-up transformation is a map from $B:= S^3 \times [-\mu,\mu]$ into $\mathbb{R}^4$. We examine the vector fields induced by this coordinate transformation in two useful coordinate charts: the entry-exit chart (or phase-directional chart) $K_1 : \{ \overline{v} = 1 \}$ and the rescaling (or central) chart $K_2 : \{ \overline{\eps} = 1 \}$, beginning with $K_1$.

In chart $K_1$, the blow-up coordinates are 
\begin{equation}  \label{K1coords}
u = r_1^2 u_1, \quad v = r_1^4, \quad \theta = r_1 \theta_1, \quad \eps = r_1^4 \eps_1,
\end{equation}
where the subscript corresponds to the chart number. The governing equations are
\begin{align} \label{fvdpK1}
\begin{split}
\dot{u}_1 &= 1 - u_1^2 - \frac{1}{3}r_1^2 u_1^3 - \frac{1}{2}\eps_1 u_1 F, \\
\dot{r}_1 &= \frac{1}{4} r_1 \eps_1 F, \\
\dot{\theta}_1  &= r_1 \eps_1 \overline{\omega} - \frac{1}{4}\eps_1 \theta_1 F, \\
\dot{\eps}_1 &= - \eps_1^2 F,
\end{split}
\end{align}
where $F(u_1,\theta_1,r_1) =-u_1 + r_1^2 \eps_1 \gamma + \beta \sqrt{\eps_1} \left( \cos(r_1\theta_1) - 1 \right)$, and we have desingularised the vector field by a factor of $r_1^2$ and recycled the overdot to denote the derivative with respect to the new time.
The hyperplanes $\{ r_1 = 0 \}$ and $\{ \eps_1 = 0 \}$ are invariant. In the invariant subspace $\{ r_1=0 \}$, $\theta_1=0$ is an attracting fixed point.
The line
\[ \ell_u = \{ (u_1, r_1, \theta_1, \eps_1)=(u_1, 0, 0, 0) \} \]
is invariant, and on it the system dynamics are governed by $\dot{u}_1 = 1 - u_1^2$. Furthermore, on $\ell_u$ there are attracting and repelling fixed points $p_a=(1,0,0,0)$ and $p_r=(-1,0,0,0)$, which, respectively, have center manifolds $N_{a,1}$ and $N_{r,1}$ in the half space $\eps_1>0$.

In order to demonstrate the existence of the primary maximal canards, we will show that there is a heteroclinic connection between $p_a$ and $p_r$ in the hyperplane $\{ r_1=0 \}$ and that this heteroclinic orbit persists for sufficiently small values of $r_1$, using Melnikov theory. The persistent connections correspond to the primary maximal canards. We carry out the relevant analysis in the rescaling chart $K_2$, where the blow-up transformation is given by
\begin{equation}  \label{K2coordsFSNI}
u = r_2^2 u_2, \quad v = r_2^4 v_2, \quad \theta = r_2 \theta_2, \quad \eps = r_2^4.
\end{equation}
Note that $r_2 = \eps^{1/4}$, so that chart $K_2$ corresponds to an $\eps$-dependent rescaling of the forced vdP equation. Also, the coordinates in the two charts are related via the following transformation:
\[ r_2 = r_1 \eps_1^{1/4}, \quad u_2 = u_1 \eps_1^{-1/2}, \quad v_2 = \eps_1^{-1}, \quad \theta_2 = \theta_1 \eps_1^{-1/4}, \]
where $\eps_1 >0$.

In chart $K_2$, the blown-up system \eqref{fvdpRegion1} is
\begin{equation} \label{fvdpK2}
\begin{split}
\dot{u}_2 &= v_2-u_2^2-\frac{1}{3}r_2^2 u_2^3, \\
\dot{v}_2 &= -u_2+r_2^2 \gamma + \beta \left(\cos(r_2\theta_2)-1\right), \\
\dot{\theta}_2 &= r_2 \, \overline{\omega},
\end{split}
\end{equation}
where we have desingularised the vector field (i.e., rescaled by $r_2^2$) and again recycled the overdot, this time to denote the derivative with respect to the new time $t_2$. This system is singularly perturbed with fast variables $(u_2,v_2)$ and slow variable $\theta_2$. Rewriting the blown-up system in non-autonomous form, we have
\begin{equation} \label{K2nonautonomous}
\begin{split}
{\dot u}_2 &= v_2-u_2^2-\frac{1}{3}r_2^2 u_2^3, \\
{\dot v}_2 &= -u_2 + r_2^2 \gamma +\beta\left(\cos(r_2^2\overline{\omega}t_2) \cos (r_2 \theta_{2,0})-1\right)-\beta \sin(r_2^2\overline{\omega}t_2) \sin(r_2\theta_{2,0}),
\end{split}
\end{equation}
where $\theta_{2,0}$ is an arbitrary phase. This is the system we will analyse in chart $K_2$. For $\overline{\omega}$ and $t_2$ of ${\cal O}(1)$, we have 
\begin{align*}
\cos(r_2^2\overline{\omega}t_2) \cos (r_2 \theta_{2,0})-1 &= \mathcal{O}(r_2^2) \text{ as } r_2 \to 0, \\
\sin(r_2^2\overline{\omega}t_2) \sin(r_2\theta_{2,0}) &= \mathcal{O}(r_2^3) \text{ as } r_2 \to 0.
\end{align*}

\begin{figure}[ht]
\centering
\includegraphics[width=3in]{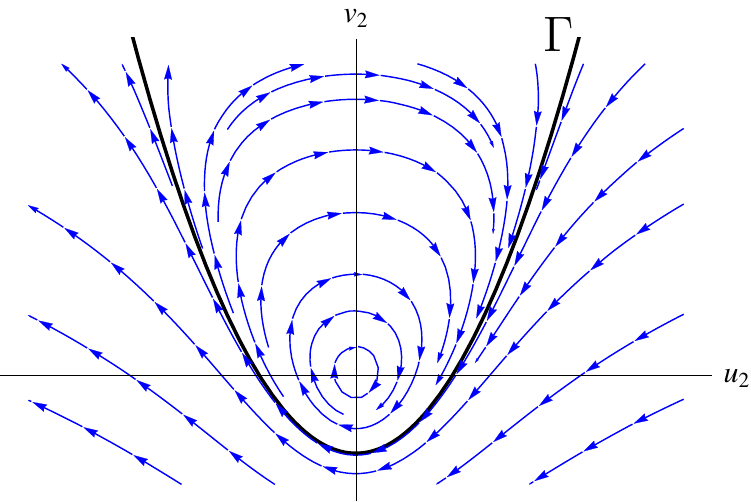}
\caption{Contours of the Hamiltonian function $H$. Periodic solutions of \eqref{unperturbedHsystem} correspond to level sets with $-\frac{1}{2}\leq H < 0$. Unbounded solutions of \eqref{unperturbedHsystem} correspond to $H>0$. The $H=0$ contour, $\Gamma$, is the heteroclinic connecting the points $p_a$ and $p_r$ at infinity.}
\label{fig:unperturbed}
\end{figure}

The unperturbed problem corresponding to \eqref{K2nonautonomous} is obtained by setting $r_2=0$,
\begin{equation}  \label{unperturbedHsystem}
\begin{split}
u^\prime_2 &= v_2-u_2^2, \\
v^\prime_2 &= -u_2.
\end{split}
\end{equation}
This system is Hamiltonian with Hamiltonian function
\[ H(u_2,v_2) = e^{-2v_2} \left( u_2^2-v_2-\frac{1}{2} \right), \]
and non-canonical formulation
\begin{align*}
\dot{u}_2 &= \frac{1}{2} e^{2v_2} \frac{\partial H}{\partial v_2}, \\
\dot{v}_2 &= -\frac{1}{2} e^{2v_2} \frac{\partial H}{\partial u_2}.
\end{align*}
The level curves of $H$ are presented in Figure~\ref{fig:unperturbed}. 
%The $H=0$ level set in particular, which we denote by $\Gamma$, separates closed trajectories from unbounded orbits and has the explicit time parametrization
The contour $\Gamma$ separates closed trajectories from unbounded orbits and has the explicit time parametrization
\[ (u_{2,\Gamma},v_{2,\Gamma}) = \left( -\frac{1}{2}t_2, \frac{1}{4}t_2^2-\frac{1}{2} \right). \]
The separatrix $\Gamma$ corresponds to the singular strong canard of the FSN I. In geometric terms, it is a heteroclinic orbit 
that lies on the upper-hemisphere and connects the fixed points $p_a$ and $p_r$, which both lie on the equator of the blown-up sphere. 

We now use the Melnikov method to analyse the persistence of $\Gamma$ under small-amplitude perturbations. As applied to \eqref{K2nonautonomous}, Melnikov theory measures the splitting distance $D$ between the curves of solutions of the perturbed system 
that are forward and backward asymptotic to $p_r$ and $p_a$, respectively. 
%For \eqref{fvdp}, the odd term proportional to sine does not contribute to the distance measurement $D$.
We develop $D$ in an asymptotic series in the small parameter $r_2$:
\[ D(r_2) = d_1 r_2^2 + d_2  r_2^3 + \cdots, \]
where the terms in the Melnikov integral are given by
\begin{align*}
d_1 &= \int_{-\infty}^{\infty} \nabla \left. H\right|_{\Gamma} \boldsymbol{\cdot} \begin{pmatrix} -\frac{1}{3}u_{2,\Gamma}^3 \\ \gamma + \frac{\beta}{r_2^2} \left( \cos(r_2^2 \overline{\omega}t_2)\cos(r_2 \theta_{2,0})-1 \right)  \end{pmatrix} \, dt_2, \\
d_2 &= \int_{-\infty}^{\infty} \nabla \left. H\right|_{\Gamma} \boldsymbol{\cdot} \begin{pmatrix} 0 \\ -\frac{\beta}{r_2^3} \sin(r_2^2 \overline{\omega}t_2)\sin(r_2\theta_{2,0}) \end{pmatrix} \, dt_2. 
\end{align*}
We note that the integrand in $d_2$ is an odd function of $t_2$ so the integral evaluates to zero and the sine term has no contribution to the distance measurement $D$.  We also note that the $\frac{\cos(r_2^2 \overline{\omega}t_2)\cos(r_2 \theta_{2,0})-1}{r_2^2}$ term in $d_1$ is $\mathcal{O}(1)$ with respect to $r_2$. The integral $d_1$ was evaluated by taking $\cos z = \operatorname{Re} (e^{iz})$, completing the square in the exponential, and deforming the contour in the complex plane. The result is
\[ d_1 = \frac{e\sqrt{2\pi}}{r_2^2} \left\{ \beta-r_2^2 \left( \frac{1}{8}+\gamma \right) - \beta e^{-\frac{1}{2} r_2^4 \overline{\omega}^2 } \cos(r_2 \theta_{2,0}) \right\}. \]
Substituting this into the bifurcation equation $D=0$, we have
\[ r_2^2 \left( \frac{1}{8} + \gamma \right) + \beta \cos(r_2 \theta_{2,0}) \left( e^{-\frac{1}{2} r_2^4 \overline{\omega}^{2}}  - 1 \right) = 0. \]

Thus, reverting to the parameters $a,b,\eps$, and $\overline{\omega}$, we see that the primary maximal canards for the FSN I are given by
\begin{equation}	\label{theta0-mainformulas-Region1}
a = 1 - \frac{\eps}{8}-b \cos(\theta_0) \exp \left(-\frac{\eps\overline{\omega}^2}{2}\right),
\end{equation}
which is \eqref{theta0-mainformulas}. We remark that $\theta_0$ is the arbitrary phase $\theta_{2,0}$ in the original $\theta$ coordinate (i.e., $\theta_0=r_2 \theta_{2,0}$). This completes the demonstration that \eqref{theta0-mainformulas} holds for low-frequency forcing, giving the locus of points at which the primary maximal canards exist. Moreover, one also sees that, for each $\overline{\omega}$,
the folds of the primary maximal canards at the endpoints of these parameter intervals are given by
\begin{equation}	\label{mainformulas-Region1}
a = 1 - \frac{\eps}{8} \pm b \exp \left(-\frac{\eps\overline{\omega}^2}{2}\right),
\end{equation}
which is precisely \eqref{mainformulas-Reg1}.
The loci in the $(\overline{\omega},a)$ plane of the folds of maximal canards mark the upper and lower boundaries of the regime in which the primary maximal canards exist.
This completes the proof of Theorem~\ref{thm:lowfrequency}.

%Finally, we note that
%the above analysis and result
%for \eqref{theta0-mainformulas} with $b={\cal O}(\sqrt{\eps})$
%can be extended from the regime in which $\omega={\cal O}(\eps)$
%to the regime in which $\omega={\cal O}(\eps^{3/4})$.
%This follows from tracking the factors of $r_2^2 \overline{\omega}$,
%and seeing that the integrals still converge for this.
%This is significant because from the numerical continuations
%we see that, in the regime where $b={\cal O}(\sqrt{\eps})$,
%the continuation effectively ends 
%when $\omega={\cal O}(\eps^{3/4})$.
%For example, with $b=0....$, $\eps=0.01$,
%we find that the continuation effectively ends
%at $\omega=...$, which is ${\cal O}(\eps^{3/4})$.

\begin{remark}
For $\omega=\eps\overline{\omega}$, one may extend the result of Theorem~\ref{thm:lowfrequency} to the parameter regime in which $b=\mathcal{O}(\eps^{\frac{1}{4}})$. Let $b=\eps^{\frac14}\beta$ and $\eta=\eps^{\frac34}\gamma$, where $\beta$ and $\gamma$ are ${\cal O}(1)$ with respect to $\eps$. Then, the perturbation terms in the $v_2$ component of the non-autonomous system become:
$$
r_2 \left[
        \gamma
        +\frac{\beta}{r_2^2}
           \left(
              \cos(r_2^2\overline{\omega}t_2)\cos(r_2\theta_{2,0})-1
           \right)
        -\frac{\beta}{r_2^2}
           \sin(r_2^2\overline{\omega}t_2)\sin(r_2\theta_{2,0})
    \right].
$$
Here, the even terms are ${\cal O}(r_2)$ as $r_2 \to 0$, so that one may
proceed with a similar Melnikov calulation as above, and the odd terms
again do not contribute to leading order in the Melnikov calculation.

We further note that a blow-up and Melnikov computation similar to that just presented for the FSN I points may be done for the folded nodes and folded saddles, and this gives the location of the maximal canards as
\[ \theta_{n,s}(k) \approx 2k \pi \pm \cos^{-1} \left( \frac{1-a-\eps/8}{b} + \mathcal{O}(b) \right), \qquad k \in \mathbb{Z}. \]
See also equation (42) in~\cite{Benes2011}.
\end{remark}

%\newpage
%---------------------------------------------------------------------------------
\section{Loci of the Torus Canards and Their Folds for Intermediate- and High-Frequency Forcing}
\label{sec:AnalysisRegion2+3}
%---------------------------------------------------------------------------------
% Calculations  of Dec 16, 2014
%---------------------------------------------------------------------------------

In this section, 
we study system \eqref{fvdp}
in the intermediate-frequency regime
with $\omega=\sqrt{\eps}\Omega$ 
and $\Omega={\cal O}(1)$,
as well as in the high-frequency regime
with $\omega={\cal O}(1)$.
We prove Theorem~\ref{thm:intermediatefrequency},
demonstrating
that system \eqref{fvdp}
possesses a family of torus canards
in the intermediate-frequency regime,
in between the two fold curves \eqref{mainformulas-Reg2}
of these torus canards.
The central methods used in the proof
are geometric desingularization
--in which we use a cylindrical blow-up of the fold curve
rather than a spherical blow-up as used
in the previous section-- and 
a Melnikov calculation 
to identify the parameter values 
for which the torus canards exist.
After Theorem~\ref{thm:intermediatefrequency} is established,
we also prove Corollary~\ref{cor:highfrequency}
for the high-frequency regime.

In the intermediate-frequency regime,
system \eqref{fvdp} is equivalent to
\begin{equation}        \label{fvdp-fasttime}
\begin{split}
x' &= y-f(x), \\
y' &= \eps(-x + a+b \cos \theta), \\
\theta' &= \sqrt{\eps} \Omega.
\end{split}
\end{equation}
First, we rectify the fold curve so that it coincides
with the $\theta$ axis,
\[ u = x-1, \quad v=y+\frac{2}{3}, \quad \alpha = a-1. \]
Also, we recall
\[
\alpha = \eps \tilde{\alpha}, \quad
b = \eps \tilde{\beta}, 
\]
where $\tilde{\alpha}$
and $\tilde{\beta}$ are ${\cal O}(1)$
with respect to $\eps$.
This transforms \eqref{fvdp-fasttime} 
to the following system:
\begin{equation}\label{u-v-theta-system}
\begin{split}
u' &= v-u^2-\frac{1}{3}u^3, \\
v' &= -\eps u + \eps^2(\tilde{\alpha} + \tilde{\beta}\cos \theta), \\
\theta' &= \sqrt{\eps} \Omega.
\end{split}
\end{equation}
Next, we perform the following cylindrical blow-up transformation:
\begin{equation}
u = \overline{r} \overline{u}, \quad 
v = \overline{r}^2 \overline{v}, \quad 
\eps = \overline{r}^2 \overline{\eps}, 
\end{equation}
which transforms the circle of fold points into a torus.
(This contrasts with the spherical blow-up
of the FSN I point in the previous section.)
Append the trivial equation
$\eps^{\prime}=0$
to system
\eqref{u-v-theta-system}
and let $\mu>0$
be sufficiently small.
For each $\theta\in \mathbb{S}^1$
and for all non-negative values 
of the system parameters,
the coordinate change is
a map from
$B := \mathbb{S}^2 \times 
[-\mu,\mu]$
into $\mathbb{R}^3$.
We examine the vector fields
induced by
\eqref{u-v-theta-system}
in two useful coordinate charts:
the entry-exit chart
(or phase-directional chart)
$K_1 = \{ \overline{v} = 1 \}$
and the rescaling chart
$K_2 = \{ \overline{\eps} = 1 \}$.

In chart $K_1$,
the coordinates are
\[
u = r_1 u_1, \quad
v = r_1^2, \quad
\eps = r_1^2 \eps_1.
\]
Setting $F(u_1, r_1, \eps_1, \theta, \tilde{\alpha}, \tilde{\beta})
= - u_1 + \eps_1 r_1 ( \tilde{\alpha} + \tilde{\beta} \cos(\theta))$,
we find that the system
in chart $K_1$ is
\begin{align}
\begin{split}
\label{eq-chartK1}
\dot{u}_1 &= 1 - u_1^2 
                 - \frac13 r_1 u_1^3
                 - \frac{u_1 \eps_1}{2} F, \\
\dot{r}_1 &= \frac{r_1 \eps_1}{2} F, \\
\dot{\theta} &= \sqrt{\eps_1}\Omega, \\
\dot{\eps}_1 &= - \eps_1^2 F,
\end{split}
\end{align}
where we have desingularised
the vector field
by rescaling the time variable
by a factor of $r_1$
and recycled the overdot.
In the phase space
of \eqref{eq-chartK1},
the hyperplanes
$\{r_1=0\}$
and $\{\eps_1=0\}$
are invariant.
In addition,
for every $(\tilde{\alpha},\tilde{\beta})$,
there is an invariant line
\[
\ell_u = \{ (u_1, r_1, \eps_1)=
(u_1,0,0) \}
\]
on which the dynamics
are governed by
$\dot{u}_1 = 1 - u_1^2$ and $\dot{\theta}=0$.
Moreover, 
for every $(\tilde{\alpha},\tilde{\beta})$,
the points
\[
p_a = ( 1, 0, 0 )
\quad
{\rm and}
\quad
p_r = ( -1, 0, 0 )
\]
are attracting and repelling fixed points,
respectively, on $\ell_u$,
and they have two-dimensional
center manifolds
$N_{a,1}$ and
$N_{r,1}$
in the half-space $\eps_1 > 0$.

In order to establish the existence
of the torus canards,
we now hook up the dynamics 
observed in chart $K_1$ 
to those in chart $K_2$.
In chart $K_2$,
the coordinates are
\begin{equation}  \label{K2coords}
u = r_2 u_2, \quad
v = r_2^2 v_2, \quad
\eps = r_2^2,
\end{equation}
and these coordinates are related to those of chart $K_1$ via the following coordinate transformation:
\[ u_2 = u_1 \eps_1^{-1/2}, \quad v_2 = \eps_1^{-1}, \quad r_2 = r_1 \eps_1^{1/2}, \]
where $\eps_1 >0$.

In chart $K_2$,
the system is
\begin{align} \label{K2-intermediatefrequency}
\dot{u}_2 &= v_2-u_2^2-\frac{1}{3}r_2 u_2^3, \nonumber \\
\dot{v}_2 &= -u_2 + r_2 \left( \tilde{\alpha}
                         + \tilde{\beta} \cos (\Omega t_2 + \theta_0) \right),
\end{align}
where we have also rescaled time by a factor of $r_2$ 
in order to desingularise the vector field 
(with $t_2$ denoting this rescaled time variable), 
recycled the overdot again now to denote the derivative 
with respect to $t_2$, 
and written the system as a non-autonomous system. 
For reference, we emphasize the relation $\sqrt{\eps}=r_2$.
We now show that system \eqref{K2-intermediatefrequency} 
possesses a special family of homoclinic orbits, 
connecting the point at infinity to itself, 
which implies that the orbits connect
the points $p_r$ and $p_a$
identified in chart $K_1$.
These orbits correspond to singular torus canards 
of the original system \eqref{fvdp-fasttime}.

The unperturbed problem associated to system 
\eqref{K2-intermediatefrequency} is
given by
\begin{align*}
\dot{u}_2 &= v_2-u_2^2, \\
\dot{v}_2 &= -u_2 ,
\end{align*}
which is the same as \eqref{unperturbedHsystem}.
As shown in the previous section,
this unperturbed system is Hamiltonian
\[ H(u_2,v_2) = e^{-2v_2} \left( u_2^2 - v_2 - \frac{1}{2} \right). \]
Along the level set $\Gamma := \{ H=0 \}$, which is the separatrix between bounded and unbounded solutions (see Figure~\ref{fig:unperturbed}), the solutions are given explicitly by
\[ u_{2,\Gamma} (t_2) = - \frac{t_2}{2}, \qquad v_{2,\Gamma} (t_2) = \frac{t_2^2}{4} - \frac{1}{2}. \]
%
%The set $\Gamma$ is a separatrix between the bounded trajectories of the unperturbed system inside and the unbounded trajectories outside. 
In the language of dynamical systems, it is a homoclinic orbit to infinity.

With the above information about the unperturbed system in hand, 
we now turn to show that $\Gamma$ persists 
for sufficiently small values of $r_2$ 
in \eqref{K2-intermediatefrequency}. 
We use a straightforward generalization 
of Proposition 3.5 of \cite{Krupa2001}, 
where we note that the perturbation terms there 
are strictly autonomous,
whereas here the perturbation terms 
also include a small-amplitude, time-periodic function,
and a compactification of the phase space can be used.
Moreover, we observe that 
the parameter there, $\lambda_2$, is also
treated as being a small variable
via the linear scaling of $\lambda$ with $r_2$,
whereas here we have chosen instead 
to scale the parameters $\alpha$ and $b$
with $\eps$ from the outset
and to treat $\tilde{\alpha},\tilde{\beta} = {\cal O}(1)$ as parameters.
In this manner, $r_2$ is the only small variable in the analysis here.

The splitting distance
between the manifolds $N_{a,2}$
and $N_{r,2}$ for system \eqref{K2-intermediatefrequency} is
\[ D(r_2) = d_{r_2}r_2 + \cdots. \]
Here, the dependence of the Melnikov function 
on the system parameters is implicit. 
We find
\begin{align}    \label{d-intermediatefrequency}
d_{r_2} &=
  \int_{-\infty}^{\infty} \nabla \left. H\right|_{\Gamma} \boldsymbol{\cdot}
          \begin{pmatrix} -\frac{1}{3}u_{2,\Gamma}^3 \\ 
                          \tilde{\alpha} + \tilde{\beta} \cos(\Omega t_2 + \theta_0)
           \end{pmatrix}
          \, dt_2 \nonumber \\
  &= -\frac{e}{2} \sqrt{2\pi} \left( \frac{1}{8} 
     + \tilde{\alpha} 
     + \tilde{\beta} e^{-\frac{\Omega^2}{2}} \cos(\theta_0) \right),
\end{align}
where the last term in the integral was evaluated 
by using $\cos(z)= {\rm Re} (e^{iz})$, 
completing the square on the exponential, 
and shifting the contour in the complex plane. 
Hence, reverting to the given parameters, 
we see that to leading order the splitting distance is
\begin{equation}        \label{D-intermediatefrequency}
D = -\frac{e}{2} \sqrt{2\pi}
\left[ \frac{\sqrt{\eps}}{8} + \frac{a - 1}{\sqrt{\eps}}
+ \frac{b}{\sqrt{\eps}} e^{-\frac{\Omega^2}{2}} \cos(\theta_0) \right].
\end{equation}

Therefore, 
for each $b$ satisfying the hypotheses
of the theorem and for $\eps$ small enough,
the simple zeroes of the Melnikov function are given 
in the $(a,\Omega)$ plane by
\begin{equation}  \label{a-Omega}
a = 1 - \frac{\eps}{8} - b e^{-\frac{\Omega^2}{2}} \cos(\theta_0).
\end{equation}
This formula, which is exactly \eqref{theta0-mainformulas}
in the intermediate-frequency regime,
gives the parameter values 
for which system \eqref{K2-intermediatefrequency} 
has a one-parameter ($\theta_0$) family 
of persistent homoclinic orbits, 
and these persistent homoclinic orbits 
of \eqref{K2-intermediatefrequency} 
are the torus canards of \eqref{fvdp-fasttime}.

Also, as a direct corollary,
we observe that the envelope of the family of torus canards
is given by
\begin{equation}  \label{sec4.2.Omega-a}
a = 1 - \frac{\eps}{8} \pm b \exp \left(-\frac{\Omega^2}{2}\right), 
\end{equation}
which is precisely formula \eqref{mainformulas-Reg2}.
This completes the proof of Theorem~\ref{thm:intermediatefrequency}.

To conclude this section, 
we prove Corollary~\ref{cor:highfrequency}.
The proof follows by extending,
in a straight forward manner, the above analysis
of the persistent homoclinics and the folds of the torus canards
in the proof of Theorem~\ref{thm:intermediatefrequency}
for the intermediate-frequency forcing regime 
to the high-frequency forcing regime.
In particular, in the high-frequency regime, 
$\omega = {\cal O}(1)$, 
which corresponds to taking 
$\Omega={\cal O}(1/\sqrt{\eps})$
in the above analysis. Following exactly along the above
calculations, we see that the geometric desingularization method yields
the same equations \eqref{K2-intermediatefrequency} in chart $K_2$, but
now the small-amplitude time-periodic forcing term 
has high-frequency
$\Omega={\cal O}(1/\sqrt{\eps})$. 
Hence, the suitable version of the
Melnikov theory is that for rapidly forced systems, and the splitting
distance along $\Gamma$ is again given by
\eqref{sec4.2.Omega-a}, which is now exponentially small in $\eps$,
since  
$\Omega={\cal O}(1/\sqrt{\eps})$. 
See for example \cite{Delshams1992}.
This completes the proof of Corollary~\ref{cor:highfrequency}.

\begin{figure}[!h]
\centering
\includegraphics[scale=0.5]{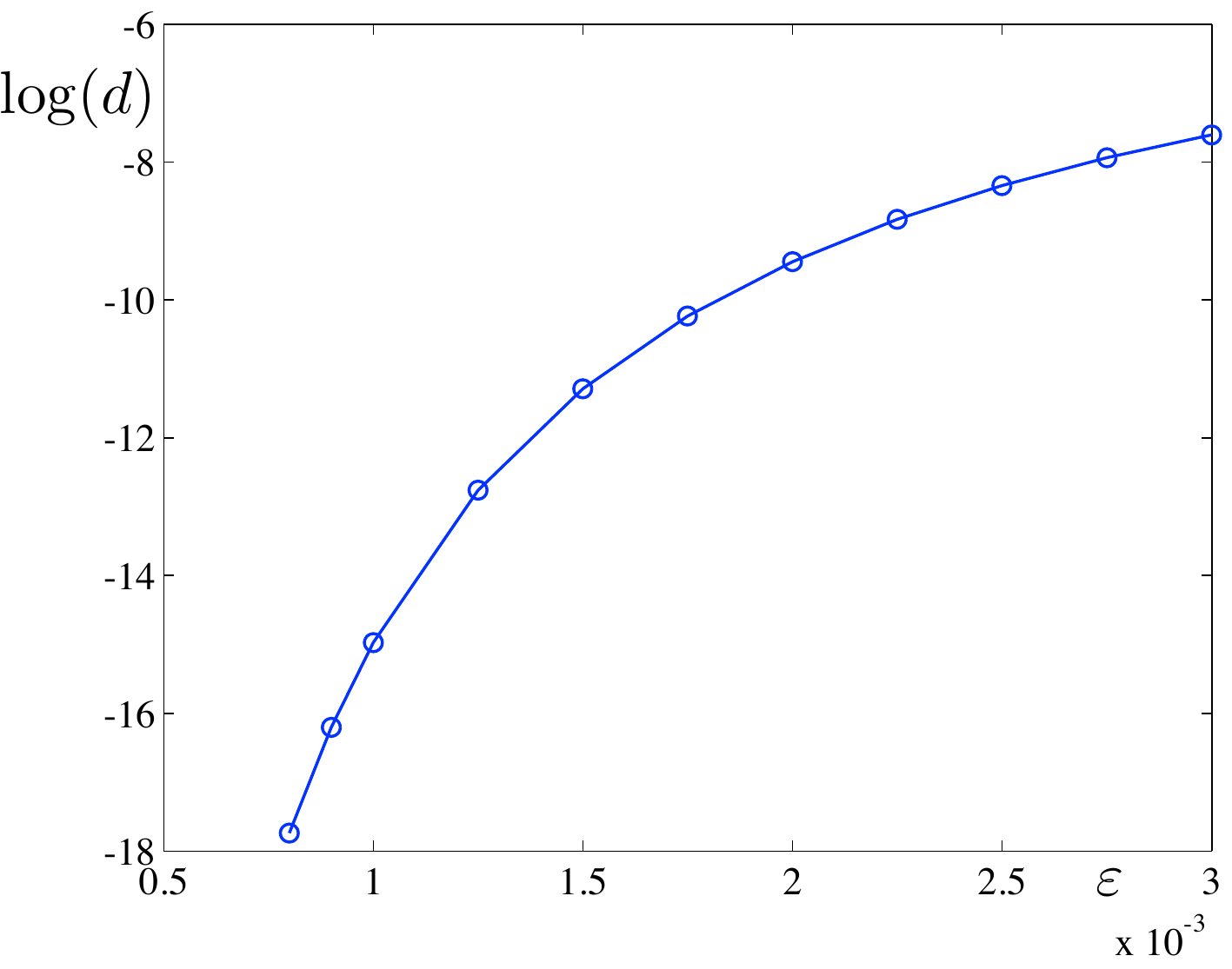}
\vspace*{-0.2cm}
\caption{\small{Width of the canard region in the $\omega=O(1)$ regime, 
as a function
of $\eps$ for $b=0.01$.}}
\label{fig:distfoldcaneps}
\end{figure}

The result of this Corollary for the high-frequency regime
also agrees well with the results obtained from
numerical simulations.
In Figure~\ref{fig:distfoldcaneps}, 
for $\omega=O(1)$,
we present a computation of the distance 
between the two folds of maximal primary canards 
as a function of $\eps$. 
We gathered the control points 
obtained for various computations 
for eleven fixed values of $\eps$, 
decreasing from $3\cdot10^{-3}$ down to $8\cdot10^{-4}$ 
and plotted them on a logarithmic scale. 
The hyperbolic shape of the resulting curve 
confrms that
this distance is exponentially small in $\eps$ as $\eps$ tends to $0$.

%Figure 5: 
%\textcolor{blue}
%{Insert Figure 5 here. This should show 
%the data points Mathieu obtained
%The latex code for including this figure
%is right below these lines, percented out.
%The label for this figure is 
%fig:distfoldcaneps 
%Also, note that this was formerly Fig 9 in the 2/19 version.}
%HERE
%
%
%THERE

%---------------------------------------------------------------------------------
\section{Secondary Canards}	\label{sec:seccan}
%---------------------------------------------------------------------------------

\begin{figure}[!h]
\centering
\includegraphics[scale=0.5]{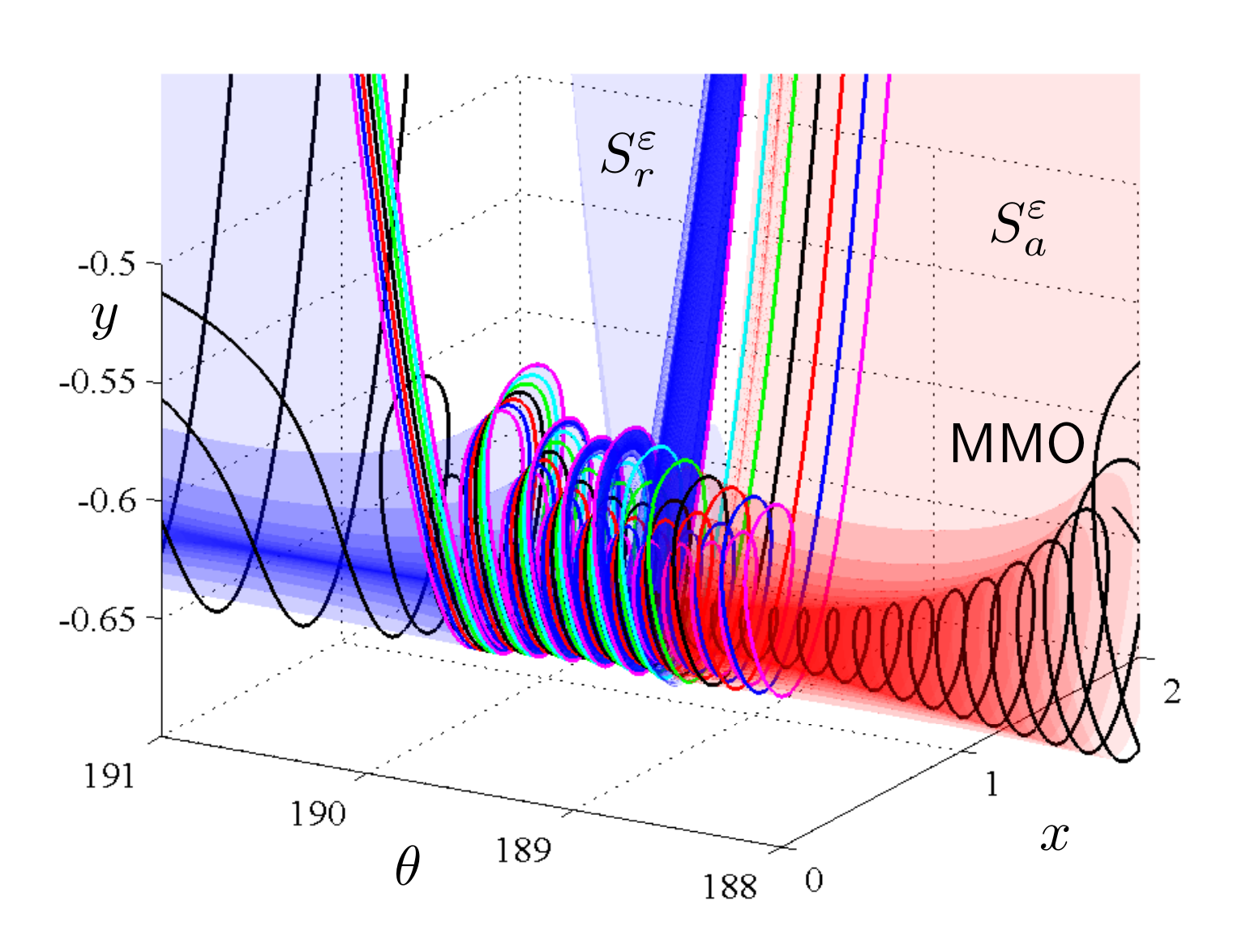}
\includegraphics[scale=0.5]{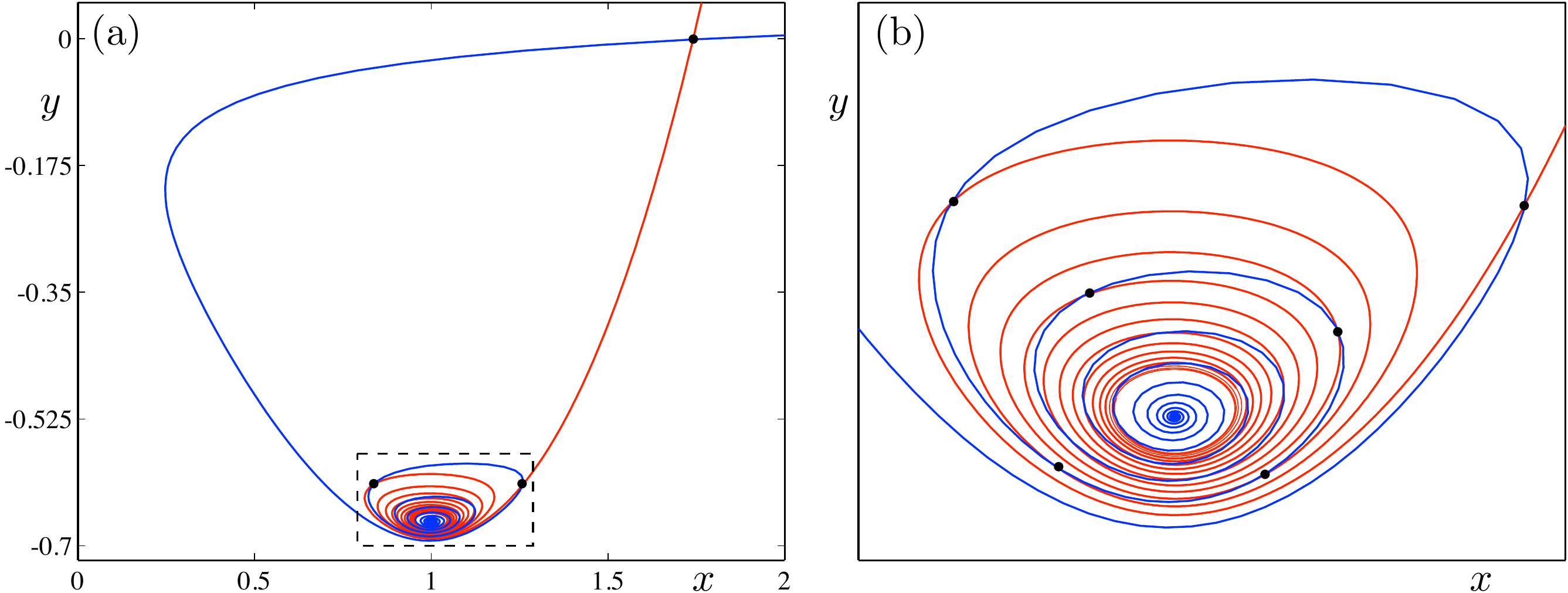}
%\vspace*{-0.2cm}
\caption{\small{Attracting (red) and repelling (blue) slow manifolds of system~\eqref{fvdp} for $a=0.9935, b=0.01, \omega=0.005$, and $\eps=0.05$, together with a stable MMO with 1 LAO and 6 SAOs found for the same parameter values by direct simulation.
Bottom row: invariant slow manifolds $S_a^{\eps}$ (red) and $S_r^{\eps}$ (blue) of system~\eqref{fvdp} in the cross-section $\Sigma_n$. The secondary canards are identified as the intersections (black dots) of $S_a^{\eps}$ and $S_r^{\eps}$.}}
\label{fig:slowman}
\end{figure}

Having established the existence of the primary strong canards and their folds, we now turn our attention to the secondary canards of the folded nodes of \eqref{fvdp}, which exist in the low-frequency forcing regime $\omega = \eps \overline{\omega}$. 
By definition, secondary canards lie in the transverse intersections of the invariant slow manifolds $S_a^{\eps}$ and $S_r^{\eps}$. A representative example of these manifolds and their intersections (i.e. the secondary canards) is shown in Figure~\ref{fig:slowman}.
These manifolds are computed from curves of initial conditions traced on the attracting and repelling sheets, respectively, of the critical manifold $S$, up to a cross-section at fixed angle $\theta$ corresponding to the maximal torus canard \cite{Desroches2008a,Desroches2010}:
\[ \Sigma_n=\left\{\theta_n=\cos^{-1}\left(\frac{1-a-\eps/8}{b}+\mathcal{O}(b)\right)\right\}. \]
%
%THERE
%
In Section \ref{subsec:resonance}, we study the folds of the secondary canards and investigate how they change in the $(\omega,a)$ plane under variation of the forcing amplitude $b$ (analogous to the folds of the primary canards). We then investigate in Section \ref{subsec:LAOtoSAO} how large-amplitude oscillations can grow from small-amplitude oscillations.

%-------------------------------------------
\subsection{Continuation of Secondary Canards}	\label{subsec:resonance}
%-------------------------------------------

As shown in Figure~\ref{fig:slowman}, the invariant slow manifolds $S_a^{\eps} \cap \Sigma_n$ and $S_r^{\eps} \cap \Sigma_n$ spiral around one another, which is typical of a folded node. 
More precisely, let $\mu := \lambda_w/\lambda_s$, $|\lambda_w|<|\lambda_s|$, denote the eigenvalue ratio of a folded node, regarded as an equilibrium of the desingularised system \eqref{desing}. Provided $\eps$ is sufficiently small and $\mu$ is bounded away from zero, the total number of (primary and secondary) maximal canards is $s_{\max}+1$, where
\[ s_{\max} = \lfloor \frac{\mu+1}{2\mu} \rfloor, \]
and $\lfloor \cdot \rfloor$ denotes the floor function. In particular, a persistent branch of secondary canards bifurcates from the weak canard in a transcritical bifurcation for odd integer values of $\mu^{-1}$ \cite{Wechselberger2005}.

\begin{remark}
The $k^{\text{th}}$ secondary canard exhibits $k$ small oscillations about the weak canard for $k=1,2,\ldots, s_{\max}-1$. These small oscillations are localized to an $\mathcal{O}(\sqrt{\eps})$ neighbourhood of the folded node \cite{Wechselberger2005,Wechselberger2012}. Moreover, trajectories on $S_a^{\eps}$ situated between $\gamma_{k-1}$ and $\gamma_{k}$, $k=1,2,\ldots,s_{\max}$ execute $k$ small oscillations about the weak canard, where $\gamma_0$ and $\gamma_{s_{\max}}$ correspond to the primary strong canard and primary weak canard, respectively.
\end{remark}

\begin{figure}[ht]
\centering
\includegraphics[width=5in]{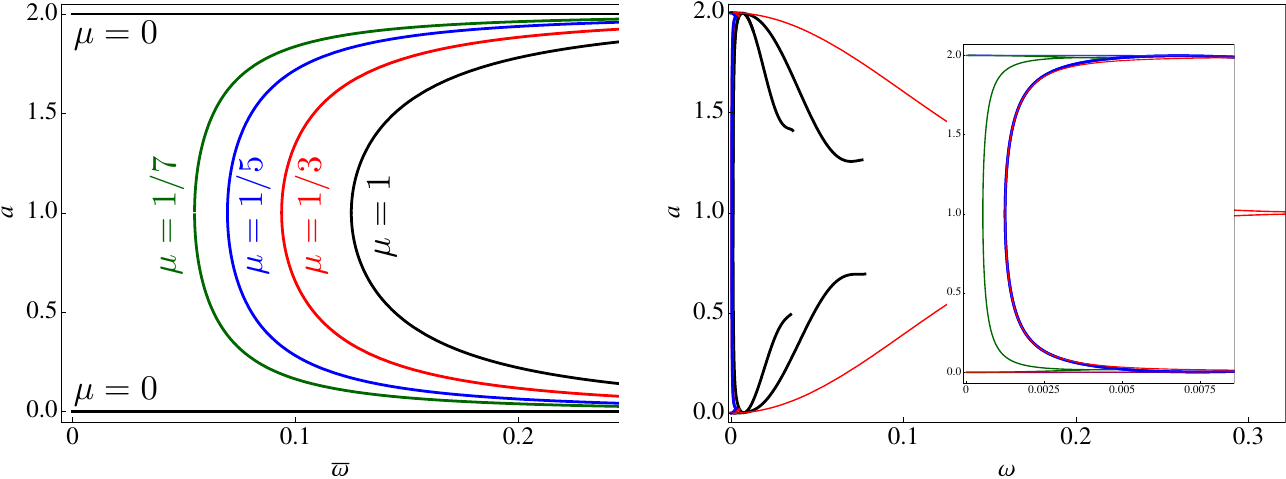}
\put(-367,126){(a)}
\put(-182,126){(b)}
\caption{\small{Resonance curves for $b=1$, and (a) $\eps=0$ and (b) $\eps=0.01$. In (b), the theoretically computed curve of maximal canards for the FSN I \eqref{mainformulas-Reg1} are shown in red. Inside this envelope, there are two black resonance curves. Both correspond to the numerical continuation of folds of maximal canards. The outermost black curve corresponds to the FSN I (i.e. where $\mu=0$). Note that \eqref{mainformulas-Reg1} breaks down when $\omega$ is no longer $\mathcal{O}(\eps)$. The inner black curve corresponds to the maximal canard of the degenerate folded node (i.e. where $\mu=1$). The inset shows the numerical continuation of the folds of canards corresponding to $\mu=1$ and $\mu=1/3$. In particular, for the $\mu=1$ resonance, we compare the numerically computed result (blue) and the theoretical result obtained in \eqref{degeneratenodeformula} (red), and find that there is excellent agreement away from the FSN I boundaries.}}
\label{fig:resonances}
\end{figure}

By tracking the resonances $\mu^{-1} = 2k+1, k=0, 1, 2, \ldots$, we can follow (in the singular limit) the locations in the $(\overline{\omega},a)$ plane where the secondary canards are born. Figure \ref{fig:resonances} shows an example for $b=1$. The non-singular $(\omega,a)$ plane shows that only the folds of canards corresponding to the FSN I and the degenerate folded node extend into the intermeditate frequency regime. All other branches of folds of canards are restricted to the low-frequency regime $\omega = \mathcal{O}(\eps)$.

\begin{remark}
Note that the resonance curves in Figure \ref{fig:resonances} bear no resemblance to the curves of folds of secondary canards in Figure \ref{fig:foldcancont}. This is to be expected since $b = \mathcal{O}(\sqrt{\eps})$ in Figure \ref{fig:foldcancont}, which implies that $\mu = \mathcal{O}(\sqrt{\eps})$ and the folded node theory does not apply.
\end{remark}

For the degenerate node ($\mu=1$), a Melnikov computation similar to that in Sections \ref{sec:AnalysisRegion1} and \ref{sec:AnalysisRegion2+3} shows that the locus of the primary maximal canard of the degenerate folded node in the $(\overline{\omega},a)$ plane is
\begin{align}\label{degeneratenodeformula}
a=1-\frac{\eps}{8}\pm \sqrt{b^2-\frac{1}{64\overline{\omega}^2}} \exp \left( -\frac{1}{2} \eps \overline{\omega}^2 \right), 
\end{align}
which holds provided $\overline{\omega}=\mathcal{O}(1), \sqrt{b^2-\frac{1}{64\overline{\omega}^2}}=\mathcal{O}(\eps)$ and $\sqrt{b^2-(1-a)^2} = \mathcal{O}(\sqrt{\eps})$. The inset of Figure \ref{fig:resonances} shows that there is excellent agreement between this theoretically computed curve and the curve obtained from numerical continuation. The deviation between theoretical and numerical results for this degenerate folded node maximal canard starts to become significant when the degenerate node branches approach the FSN I branches. We note an important implication: all secondary canards due to folded nodes are restricted to the region of the $(\omega,a)$ plane bounded by $\omega=0$, the locus of the folds of maximal canards of the FSN I and the locus of the maximal canards of the degenerate folded node. 

As was the case for the numerical continuation of the folds of the primary canards and the folds of the torus canards, the numerical continuation of the maximal secondary canards is done by solving families of boundary-value problems and computing branches of such solutions using pseudo-arclength continuation. Along these branches, a number of fold points can be detected,
and then the curves of folds of secondary canards 
can be continued in two parameters. 
For a representative set of parameter values,
the folds of the first, second, $\ldots$, tenth secondary canards (i.e., with respectively one, two, $\ldots$, ten loops) are shown in the $(\omega,a)$ plane in Figure~\ref{fig:foldcancont}. The outermost envelope in Figure~\ref{fig:foldcancont} are the curves of folds of primary canards. The rightmost path enclosed by the fold curve of the primary strong canards represents the fold curve of the first (1-loop) secondary canard, and each successive curve to the left represents a family of folds of secondary canards with one additional loop. These branches of folds of secondary canards emanate from the FSN I points at $a= 1-\frac{\eps}{8} \pm b$. As $\omega$ is increased, the corresponding pairs of $n$-loop branches come together at turning points.

\begin{figure}[!h]
\centering
\includegraphics[scale=0.5]{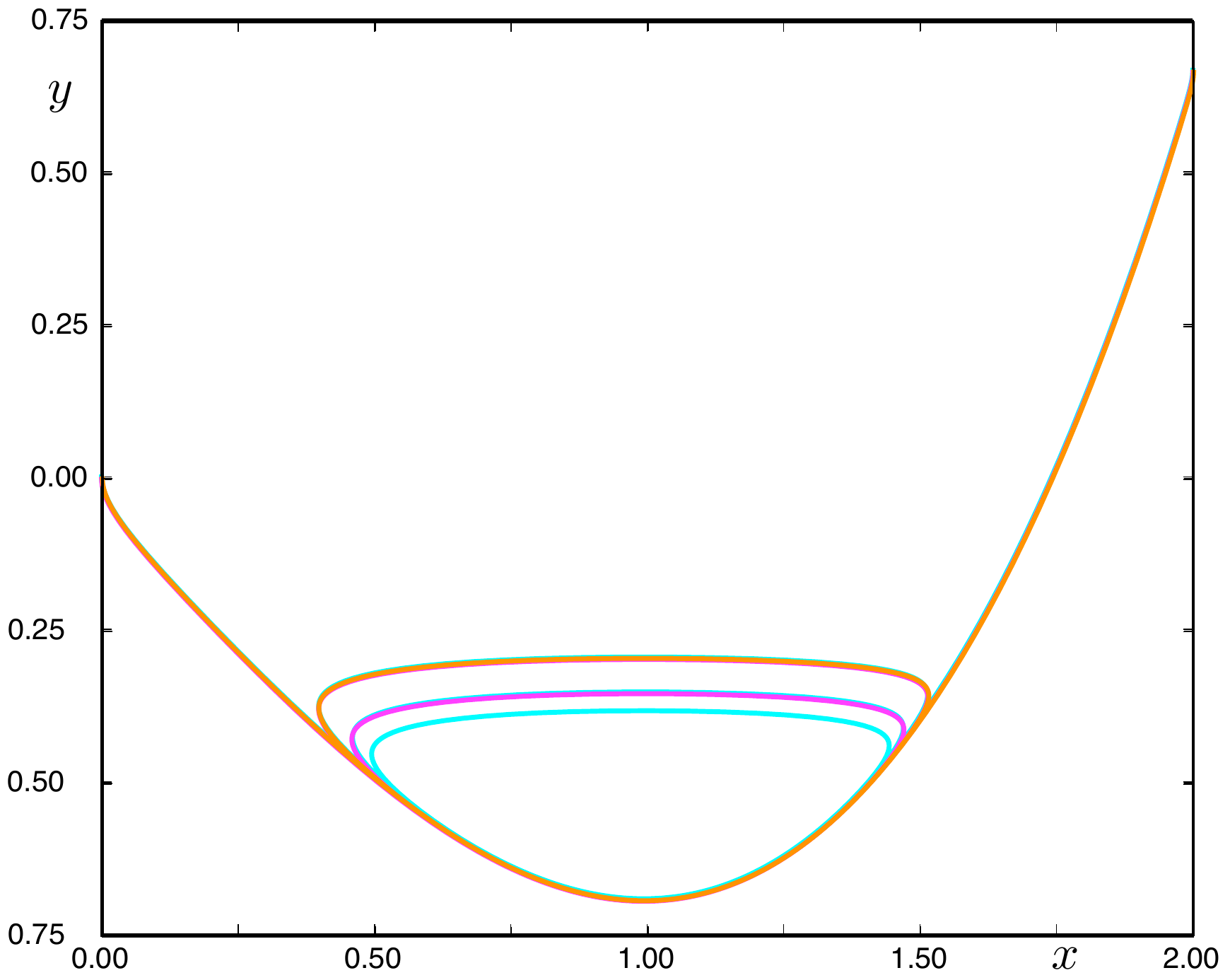}
\vspace*{-0.25cm}
\caption{\small{
Projection of the 1-loop (yellow), 2-loop (red), and 3-loop (blue) secondary canards of \eqref{fvdp} shown for $\omega=0.01$, $\eps=0.01$, $b=0.01$, and $a=1.0034909031$ (1-loop), $a=1.0034909029$ (2-loop), and $a=1.0034909029$ (3-loop) shown in the $(x,y)$-plane.}}
\label{fig:3secondarycanards}
\end{figure}

It is also useful to examine projections of the secondary canards onto the $(x,y)$ plane. In Figure~\ref{fig:3secondarycanards}, 
we show the first three maximal secondary canards, with respectively, one loop (yellow), two loops (red), and three loops (blue). The highest loops of the 2-loop and 3-loop maximal secondary canards are observed to lie extremely close to the single loop of the 1-loop maximal secondary canard. The same holds for all of the higher-loop secondary canards, as well.
Also, the second loop of the 2-loop canard
lies inside the first loop, 
and it lies extremely close to the second loop of the 3-loop canard.
In addition, for even smaller values of the forcing frequency $\omega$,
the $y$-intercepts of the return jumps increase.
Moreover, these
$y$-intercepts diverge to $\infty$ in the limit $\omega\to 0$. In fact, in this limit, the maximal secondary canards collapse onto the primary strong canard, consistent with the observation that the branches of maximal secondary canards emanate from the same FSN I points as the primary canards do.

%Figure NEW (between 6 and 7): 
%\textcolor{blue}
%{Insert a new Figure here. This should show 
% Mathieu's calculations of the 1-loop, 2-loop, and 3-loop secondary canards. 
%The label for this figure is 
%fig:3secondarycanards
%This is the new figure from March 14.}
%HERE
%
%
%THERE

\begin{figure}[!h]
\centering
\includegraphics[scale=0.5]{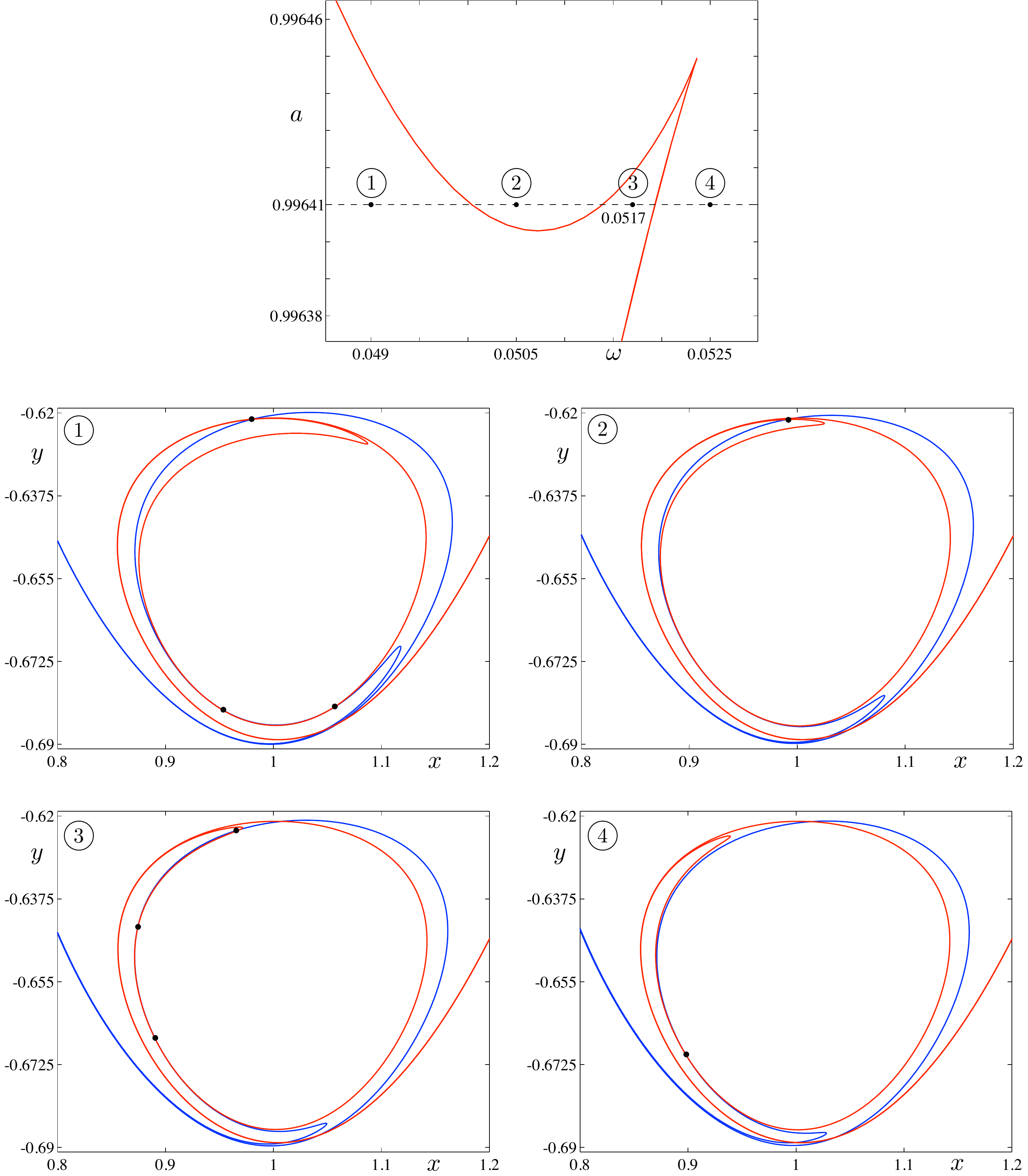}
\vspace*{-0.2cm}
\caption{\small{Evolution of the intersection points (black dots) between the curves
representing $S_a^{\eps}$ and $S_r^{\eps}$ in a fixed cross-section through
the folded node. Here, $a=0.99641$, $b=0.01$, and $\eps=0.05$. The values
of $\omega$ are (1) 0.049, (2) 0.0505, (3) 0.0517, (4) 0.0525, as also
labeled on the horizontal axis in the top frame. In the transition between
frames (1) and (2), the lower two intersection points disappear. Between
frames (2) and (3), two intersection points are created. Then, two
disappear in the final transition shown, from (3) to (4)}}
\label{fig:transcusp}
\end{figure}

We now investigate what happens to the invariant slow manifolds near the turning points (recall Figure~\ref{fig:foldcancont})
of the fold curves of secondary canards. In Figure~\ref{fig:transcusp}, we show one such fold curve and take four values of $\omega$, 
for a fixed value of $a$, near the turning point which marks the largest $\omega$-value of this curve (top panel). For each value of $\omega$, we compute $S_a^{\eps}$ and $S_r^{\eps}$ up to a fixed cross-section, 
following the procedure described above. Then, the intersection curves of both manifolds in the fixed cross-section are shown for each $\omega$-value in the four bottom panels of Figure~\ref{fig:transcusp}. Each time the fold curve of maximal canard solutions is crossed, 
two intersections of the attracting and repelling slow manifolds disappear or are created. This is illustrated in each of the transitions shown in panels (1)-(4).

%Figure 9: 
%\textcolor{blue}
%{Insert Figure 9 here. This should show 
%The latex code for including this figure
%is right below these lines, percented out.
%The label for this figure is 
%fig:transcusp
%This was formerly Fig 8 in the 2/19 version.}
%HERE
%
%
%THERE

%-------------------------------------------
\subsection{Growth of LAOs from SAOs in the Secondary Canards}	\label{subsec:LAOtoSAO}
%-------------------------------------------

\begin{figure}[!h]
\centering
\includegraphics[scale=0.5]{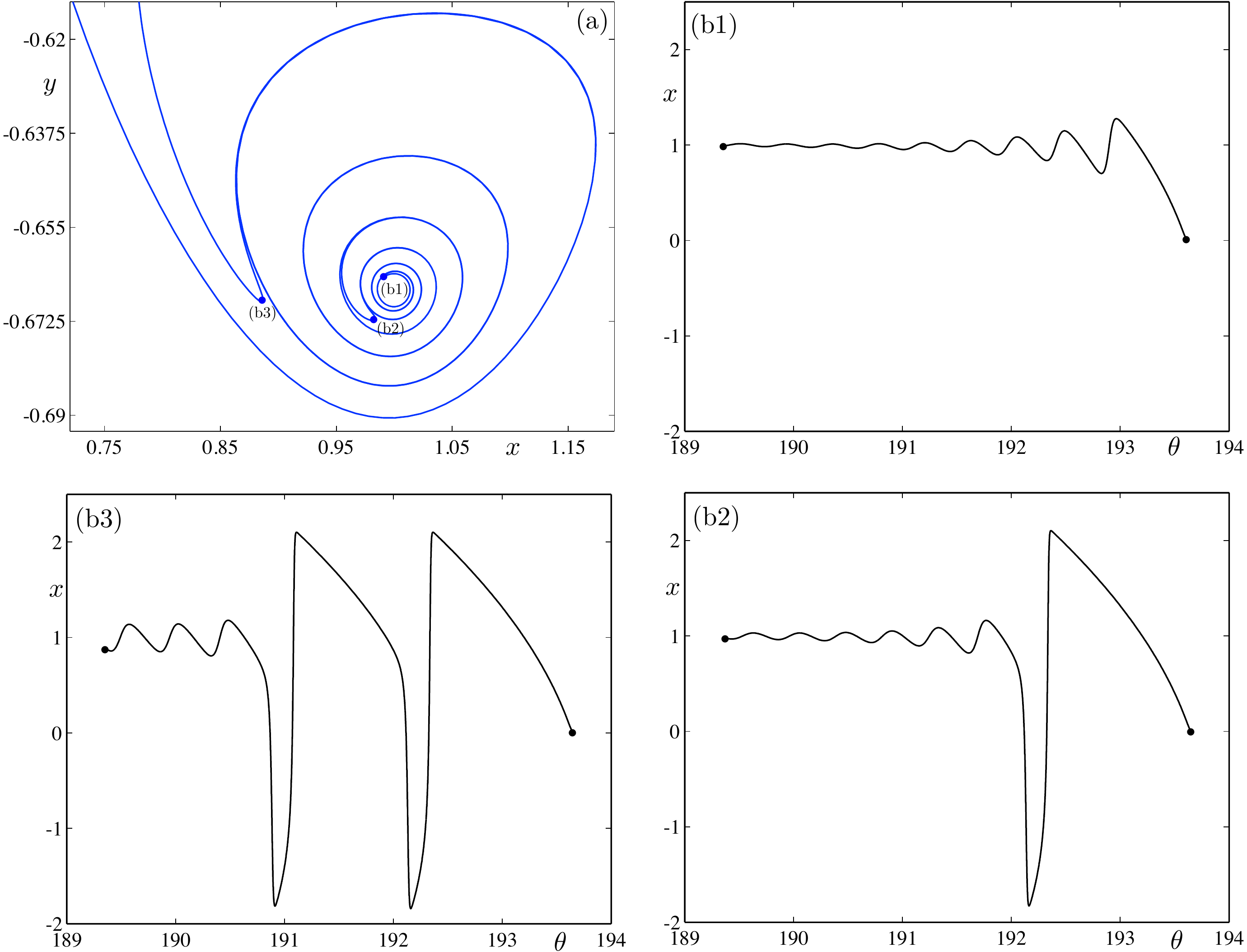}
\vspace*{-0.2cm}
\caption{\small{Repelling slow manifold reversing the direction of spiralling as successive LAOs appear in the orbit. Panels (b1), (b2), and (b3) show the solution profiles corresponding to three such events. Here, $a=0.9935, b=0.01, \omega=0.3$, and $\eps=0.05$.}}
\label{fig:LAOgrow}
\end{figure}

Along the continuation of the secondary canards, an orbit segment can `grow' an LAO. This occurs in regions where the repelling slow manifold spirals backwards instead continuing to spiral inwards towards the weak canard. In Figure~\ref{fig:LAOgrow}, we show the repelling slow manifold for $\overline{\omega}=0.3$ in the cross-section $\Sigma_n$; the direction of spiralling changes three times along this portion of the slow manifold, at the points labeled (b1), (b2), and (b3) in frame (a). Each direction reversal corresponds to the orbit growing a large-amplitude oscillation. In panels (b1), (b2), and (b3), we show the profiles of the computed solution segments at these events. Note that the first fold encountered (starting from the center of the spiral and going outwards) corresponds to the orbit segment having one SAO grow up to the size of an LAO (this occurs at the second fold), while the other SAOs stop growing in size. 

%Figure 8: 
%\textcolor{blue}
%{Insert Figure 8 here. This should show 
%evolution of intersection curves as $\omega$ is varied
%across the secondary canard points.
%The latex code for including this figure
%is right below these lines, percented out.
%The label for this figure is 
%fig:LAOgrow
%This was formerly Fig 7 in the 2/19 version.}
%HERE
%
%
%THERE

\begin{remark}
An important consequence of studying the curves of maximal canards and maximal torus canards in the parameter space of \eqref{fvdp} is that they serve as the boundaries between different dynamic regimes of \eqref{fvdp}. As highlighted briefly by the graphical summary in the $(\omega,a)$ plane shown in Figure \ref{fig:motivation}, the fvdP equation \eqref{fvdp} exhibits small-amplitude oscillations (SAOs), large-amplitude or relaxation oscillations (LAOs), and mixed-mode oscillations (MMOs). The SAOs are the $\frac{2\pi}{\omega}$-periodic solutions generated when an attracting equilibrium of the unforced vdP equation (i.e., $b=0$) is subjected to a small-amplitude periodic forcing of frequency $\omega$ (Figure \ref{fig:mmos}(a)). The LAOs occur when the equilibrium of the planar vdP equation sits on the middle branch of the cubic-shaped nullcline and the attractor of the system is a relaxation oscillation that alternates the trajectory between the outer branches of the cubic (Figure \ref{fig:mmos}(f)). The MMOs feature SAOs superimposed on large-amplitude relaxation-type oscillations. Figure \ref{fig:motivation} shows that the MMOs become more robust for low-frequency forcing, just as was observed in numerical simulations of a rotated van der Pol-type model in \cite{Benes2011}. 

\begin{figure}[h]
\centering
\includegraphics[width=5in]{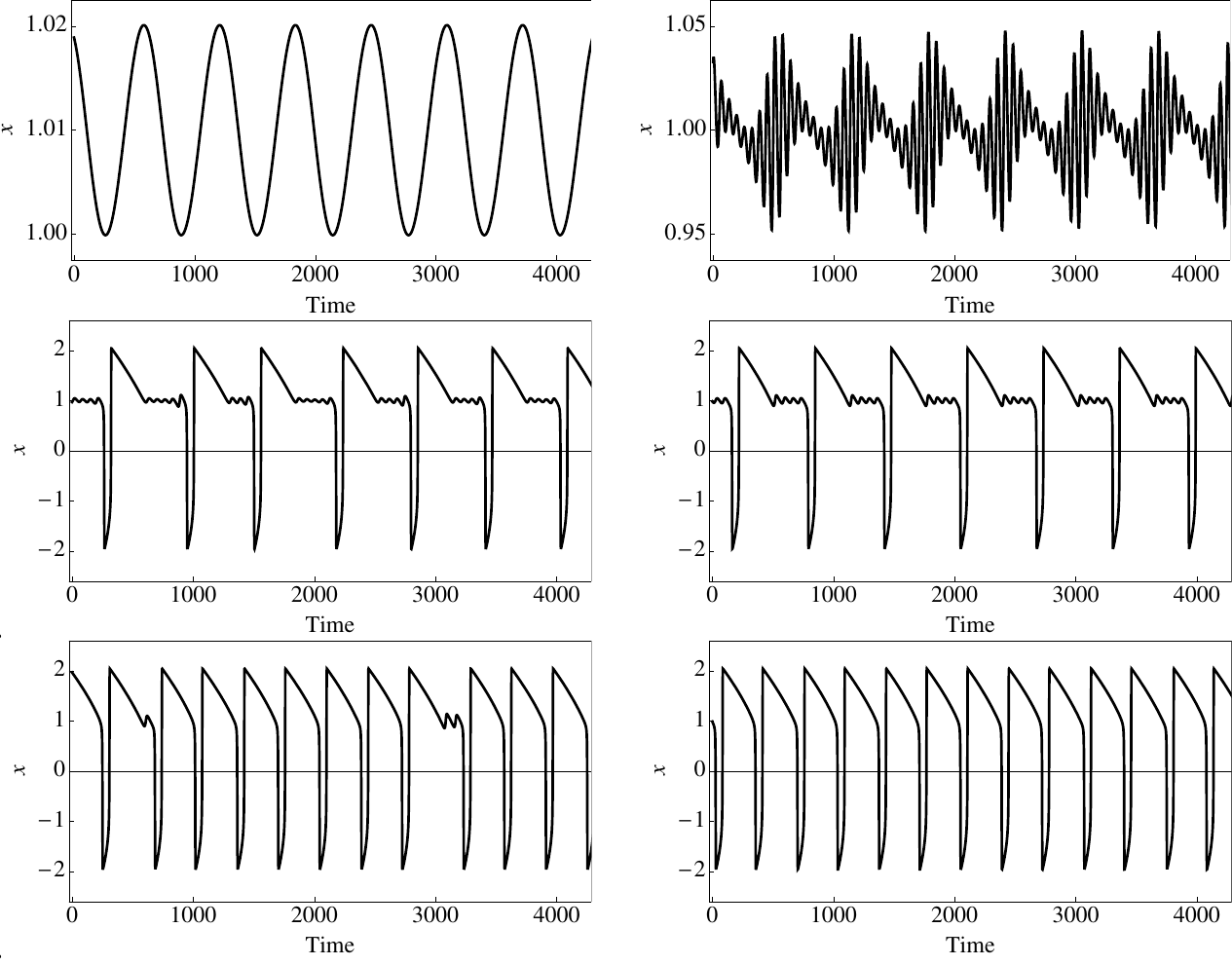}
\put(-364,274){(a)}
\put(-182,274){(b)}
\put(-364,180){(c)}
\put(-182,180){(d)}
\put(-364,84){(e)}
\put(-182,84){(f)}
\caption{Transition from SAOs to MMOs to LAOs in the forced vdP equation \eqref{fvdp}. The attractors of \eqref{fvdp} are shown for $\eps=0.01$, $\omega=0.01$, and (a) $a=1.01$, (b) $a=0.9999$, (c) $a=0.999$, (d) $a=0.994$, (e) $a=0.991$, and (f) $a=0.987$.}
\label{fig:mmos}
\end{figure}
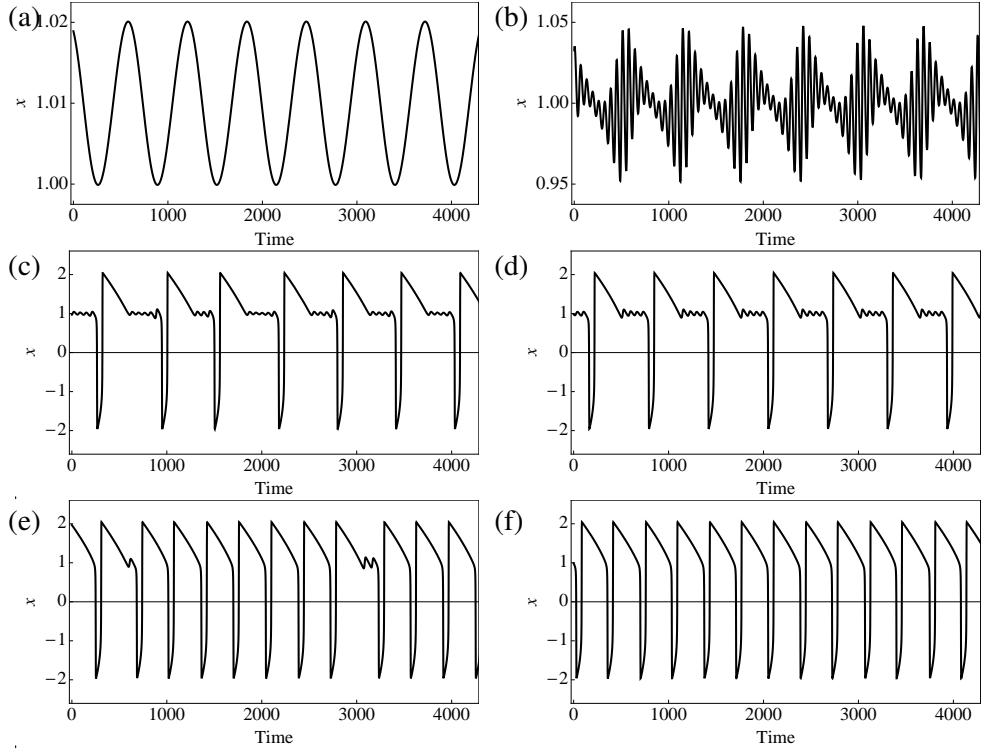

\end{remark}

%---------------------------------------------------------------------------------
\section{The fvdP is a Normal Form for a Class of Systems Near Torus Canard Explosions}
\label{sec:normalform}
%---------------------------------------------------------------------------------

In this section, we prove that, under a number of natural conditions, a slow/fast system with two fast variables and one slow variable, which is subject to time-periodic forcing and for which the fast system possesses a generic fold of limit cycles, is equivalent 
to a system in which the fast component is
given to lowest order by the forced van der Pol system \eqref{fvdp}. 
We consider systems with two fast variables and one slow variable
of the following form:
\begin{align}\label{eq-osys}
\begin{split}
\dot x &= f_1(x,y,z)\\
\dot y &= f_2(x,y,z)\\
\dot z &= \eps g(x,y,z),\qquad x,y,z \in \mathbb{R}.
\end{split}
\end{align}
The fast system is
\begin{align}\label{eq-osysf}
\begin{split}
\dot x &= f_1(x,y,z)\\
\dot y &= f_2(x,y,z)
\end{split}
\end{align}
in which $z$ is a parameter,
and we make the following hypotheses about the fast system:

\noindent
{\bf (H1)} 
For $z=0$, there exists a non-degenerate periodic solution ${\tilde \Gamma}$; and, here
by non-degenerate, we mean a periodic solution with finite period.

\noindent
{\bf (H2)} The Floquet multiplier of this periodic orbit at $z=0$ is one.

We prove the following theorem:
\begin{theorem}\label{th-red}
Given (H1), (H2), and a non-degeneracy assumption (see 
\eqref{eq-addhyp} below),
system \eqref{eq-osys} is locally 
(in a small neighborhood of ${\tilde \Gamma}$) 
orbitally equivalent to
\begin{align}\label{eq-vdpl}
\begin{split}
\dot \rho &= z-\rho^2+O(\rho^3)+O(\eps)\\
\dot  \theta &= 1\\
\dot z&= \eps \tilde g(\rho,\theta,z),\qquad \rho,z,\theta \in \mathbb{R},
\end{split}
\end{align}
where $\tilde g$ is $2\pi$-periodic in $\theta$,
and where $z$ and $\rho$ are small.
\end{theorem}
The fast subsystem of \eqref{eq-vdpl} 
has very similar dynamics to that of \eqref{fvdp}, 
however the slow system of the general systems
may be much richer than the slow system of \eqref{fvdp}.
We make
a more detailed comparison 
between the full systems later in this section. 
We first analyze the fast systems.

We begin the proof of Theorem \ref{th-red} 
with some preliminary transformations.
On the basis of (H1), 
one may rectify the flow of \eqref{eq-osysf} 
so that the periodic orbit becomes the unit circle.
Next, using the coordinates 
\[
x=(1+r)\cos(\theta) \quad {\rm and} \quad  y=(1+r)\sin(\theta),
\]
one may transform \eqref{eq-osysf} to 
\begin{align}\label{eq-osysfp}
\begin{split}
\dot r&= \tilde f_1(r,\theta,z)\\
\dot \theta &= \tilde f_2(r,\theta,z),
\end{split}
\end{align}
where $\tilde f_i$ is $2\pi$ periodic in $\theta$, 
$\tilde f_1(0,\theta, 0)=0$,
and $\tilde f_2\neq 0$ 
in a neighborhood of $(0,\theta, 0)$. 
Also, one may scale the time variable
so that \eqref{eq-osysfp} becomes
\begin{align}\label{eq-osysfps}
\begin{split}
\dot r&= F(r,\theta,z)\\
\dot \theta &= 1,
\end{split}
\end{align}
with $F=\tilde f_1/\tilde f_2$. 
This is a useful formulation of the fast system,
and we work directly with this throughout the proof.

In order to analyse the dynamics 
of this system for small values of $r$,
we expand:
\begin{equation}\label{eq-expaF}
F(r,\theta,z)=\psi_0(\theta,z)+r\psi_1(\theta,z)+r^2\psi_2(\theta,z)+\ldots +r^N\psi_N(\theta,z)+O(r^{N+1}).
\end{equation}
Condition (H2) is now equivalent to
\begin{equation}\label{eq-H2mod}
\int_0^{2\pi} \psi_1(\theta,0)d\theta=0,
\end{equation}
because
\[
\frac{\partial\bar F}{\partial r}(0,0)=0,
\]
where the average of $F$ is defined as
\[
\bar F(r,z)=\frac{1}{2\pi}\int_0^{2\pi} F(r,\theta, z)d\theta.
\]
Also, to analyze the dynamics for small values of $r$,
it is useful to introduce 
a new variable $\rho$ by 
\begin{equation}\label{eq-chaco}
r=\phi_0+\rho e^{\phi_1}+\rho^2\phi_2+\ldots +\rho^N\phi_N,\quad N\ge 2.
%\phi(\theta)=\int_0^\theta \psi(\nu)d\nu.
\end{equation}
where $\phi_i$ are functions of $(\theta,z)$ to be determined. 

The following lemma 
lies at the heart of the proof of this theorem:
\begin{lemma}\label{lem-chaco}
There exists a choice of 
functions $\phi_i$, $i=1,\ldots , N$,
in the coordinate change \eqref{eq-chaco}
and a function $G(\rho, z)$ such that, 
in the coordinates $(\rho, \theta)$, 
the fast system \eqref{eq-osysfps} has the form
\begin{align}\label{eq-aosysfps}
\begin{split}
\dot \rho &= G(\rho,z)+O(\rho^{N+1})\\
\dot \theta &= 1.
\end{split}
\end{align}
\end{lemma}

\begin{remark}
Hypothesis (H2) is not needed for this lemma.
\end{remark}

\noindent
{\bf Proof.} The proof is split into three steps.
First, we consider the simplest case in which $N=1$
in \eqref{eq-chaco}.
Next, we prove the Theorem for $N=2$, and 
finally we prove it for general $N$
in \eqref{eq-chaco}.

\vskip0.3truein

\noindent
{\bf Step 1.}
With $N=1$, 
the relevant transformation \eqref{eq-chaco} between
$r$ and $\rho$ is:
\begin{equation}\label{eq-chacosi}
r=\phi_0+\rho e^{\phi_1}.
\end{equation}
Differentiating \eqref{eq-chacosi} 
with respect to time,
substituting in
\eqref{eq-osysfps}, 
and Taylor expanding $F(\phi_0+\rho e^{\phi_1}, \theta, z)$
about $\phi_0$, we obtain
\begin{equation}\label{eq-step1}
\dot\rho=
(-\phi_{0,\theta}+F(\phi_0,\theta,z))
e^{-\phi_1}
+\rho(-\phi_{1,\theta}+F_r(\phi_0,\theta,z))+O(\rho^2).
\end{equation}
Hence, based on the form of the two terms in parentheses
in the right member of this equation,
we are naturally led to study the following 
system of two ordinary differential equations
with two unknown parameters:
\begin{equation}\label{eq-ODE}
\begin{split}
\phi_{0,\theta}&=F(\phi_0,\theta,z)+\lambda_0e^{\phi_1}\\
\phi_{1,\theta}&=F_r(\phi_0,\theta,z)+\lambda_1
\end{split}
\end{equation}
In this manner, proving the result for $N=1$ 
is now equivalent 
to finding $\lambda_0$ and $\lambda_1$ 
for which there exist $\phi_0$ and $\phi_1$ 
which are $2\pi$ periodic and 
which satisfy the equations \eqref{eq-ODE}. 

According to general theory of differential equations,
for every pair $(\lambda_0,\lambda_1)$,
there exists a unique solution $(\phi_0,\phi_1)$ of \eqref{eq-ODE} 
satisfying the initial conditions $\phi_0(0)=0$ and $\phi_1(0)=0$. 
Such solutions must satisfy the integral equations,
\begin{equation}\label{eq-integeq}
\begin{split}
\phi_{0}(\theta)&=\int_0^\theta F(\phi_0(\nu),\nu,z)d\nu+\lambda_0\int_0^{\theta} e^{\phi_1(\nu)}d\nu\\
\phi_{1}(\theta)&=\int_0^\theta F_r(\phi_0(\nu),\nu,z)d\nu +\lambda_1\theta.
\end{split}
\end{equation}
Also, in terms of these integral equations,
the conditions for periodicity are
\begin{equation}\label{eq-integeq-periodicity}
\begin{split}
0&=\int_0^{2\pi}F(\phi_0(\nu),\nu,z)d\nu+\lambda_0\int_0^{2\pi} e^{\phi_1(\nu)}d\nu\\
0&=\int_0^{2\pi} F_r(\phi_0(\nu),\nu,z)d\nu +2\pi\lambda_1.
\end{split}
\end{equation}

Based on the above formulation
of the integral equations, 
it is useful to define ${\cal H}\; :\;\mathbb{R}\times\mathbb{R}^2\to \mathbb{R}^2$,
as follows:
\begin{equation}   \label{eq-HN1}
{\cal H}(z,\lambda_0,\lambda_1)=\left (\begin{array}{c} {\cal H}_0(z,\lambda_0,\lambda_1)\\{\cal H}_1(z,\lambda_0,\lambda_1)\end{array}\right )=\left (\begin{array}{c}\int_0^{2\pi}F(\phi_0(\nu),\nu,z)d\nu+\lambda_0\int_0^{2\pi} e^{\phi_1(\nu)}d\nu\\
\int_0^{2\pi} F_r(\phi_0(\nu),\nu,z)d\nu +2\pi\lambda_1\end{array}\right ).
\end{equation}
Now, by the assumptions (H1) and (H2) (see \eqref{eq-H2mod}), 
$\lambda_0=0$, $\lambda_1=0$, and $z=0$ 
is a solution of \eqref{eq-integeq}
with $\phi_0=0$ and $\phi_{1,0}=\int_0^\theta F_r(0,\nu,0)d\nu$.
We now verify that the assumptions of the Implicit Function Theorem
are satisfied to show that there is a branch of nontrivial solutions
emanating from this trivial solution.
In particular, we verify that $D{\cal H}(0,0,0)$ is non-singular,
by showing that ${\rm det} D{\cal H} (0,0,0) = (2\pi)^2$.

From the definition,
we see that
${\cal H}_{1,\lambda_1}(0,0,0)=2\pi$.
We will now show that 
${\cal H}_{0,\lambda_1}(0,0,0)=0$
and ${\cal H}_{0,\lambda_0}(0,0,0)=2\pi$.
To show that ${\cal H}_{0,\lambda_1}(0,0,0)=0$, 
we start with
\[
{\cal H}_{0,\lambda_1}(0,0,0)=\int_0^{2\pi}F_r(0,\theta,0)\phi_{0,\lambda_1}(\theta)d\theta.
\]
Then, using \eqref{eq-ODE}, we obtain
\[
{\cal H}_{0,\lambda_1}(0,0,0)=\int_0^{2\pi}\frac{d}{d\theta} \phi_{0,\lambda_1}(\theta)d\theta=\phi_{0,\lambda_1}(2\pi)-\phi_{0,\lambda_1}(0).
\]
By assumption (above formula \eqref{eq-integeq}),
$\phi_0(0)=0$ and $\phi_1(0)=0$ 
for $\lambda_0$ and $\lambda_1$. 
Hence, $\phi_{0,\lambda_j}(0)=0$ for $j=1,2$.  
Therefore, we conclude from \eqref{eq-ODE} 
that $\phi_{0,\lambda_1}\equiv 0$,
so that this entry of the Jacobian of ${\cal H}$ vanishes,
as claimed.

Next, we show that ${\cal H}_{0,\lambda_0}(0,0,0) = 2\pi$. 
By an argument similar to the above,
\[
{\cal H}_{0,\lambda_0}(0,0,0)=\phi_{0,\lambda_0}(2\pi)-\phi_{0,\lambda_0}(0)=\phi_{0,\lambda_0}(2\pi).
\]
Therefore, from \eqref{eq-ODE}, we obtain
\begin{equation}\label{calcul1}
\frac{d}{d\theta}\left (\phi_{0,\lambda_0}\right )=F_r(0,\theta,0)\phi_{0,\lambda_0}+e^{\phi_1}.
\end{equation}
Observe that, for $\lambda_0=0$ and $\lambda_1=0$,
we have $d\phi_1/d\theta= F_r(0,\theta,0)$. 
Hence, \eqref{calcul1} is equivalent to
\begin{equation}\label{calcul2}
\frac{d}{d\theta}\left (\phi_{0,\lambda_0}\right )
=\frac{d \phi_1}{d\theta}\phi_{0,\lambda_0}+e^{\phi_1},
\end{equation}
from which it follows that
\[
\frac{d}{d\theta}\left (\phi_{0,\lambda_0}e^{-\phi_1}\right )=1.
\]
Finally, since $\phi_{0, \lambda_0}(0)=0$,
it follows from \eqref{eq-H2mod} 
that $\phi_{0, \lambda_0}(2\pi)=2\pi$. 
Therefore,
\[
\det(D{\cal H}(0,0))=(2\pi)^2,
\]
and by the Implicit Function Theorem
there is a branch of periodic solutions $\phi_0$ and $\phi_1$
for each $(\lambda_0, \lambda_1)$ sufficiently small,
emanating from the trivial solution.
This completes the proof of the Lemma
for the case $N=1$.

\vskip0.3truein
\noindent
{\bf Step 2.}
We now show that the lemma holds for $N=2$ in \eqref{eq-chaco}.
All quantities are expanded up to and including $\rho^2$.
For the vector field $F$, we have
\begin{equation}\label{eq-expan2}
F(r,\theta,z)=F(\phi_0,\theta,z)+\rho F_r(\phi_0,\theta,z)e^{\phi_1}+\rho^2(F_r(\phi_0,\theta,z)\phi_2+\frac12F_{rr}(\phi_0,\theta,z)e^{2\phi_1}).
\end{equation}
Also, differentiating \eqref{eq-chaco} for $N=2$ with respect to $t$,
we find
\begin{equation}\label{eq-rdot}
\dot r=\dot\rho(e^{\phi_1}+2\rho\phi_2)+\phi_{0,\theta}+\rho e^{\phi_1}\phi_{1,\theta}+\rho^2\phi_{2,\theta}.
\end{equation}
Hence, combining \eqref{eq-expan2} and \eqref{eq-rdot}, we get
\begin{equation}\label{eq-sysn2}
\begin{split}
\dot\rho(1+2\rho\phi_2e^{-\phi_1})=&(-\phi_{0,\theta}+F(\phi_0,\theta,z))e^{-\phi_1}+\rho(-\phi_{1,\theta}+ F_r(\phi_0,\theta,z))\\&+\rho^2(-\phi_{2,\theta}+F_r(\phi_0,\theta,z)\phi_2+\frac12F_{rr}(\phi_0,\theta,z)e^{2\phi_1})e^{-\phi_1}.
\end{split}
\end{equation}

Now, after multiplying both sides of \eqref{eq-sysn2}
by $(1 + 2 \rho \phi_2 e^{-\phi_1})^{-1}$,
we examine the structure of the terms
at each order of $\rho^0, \rho^1,$ and $\rho^2$.
This suggests that we analyze 
the following system of differential equations:
\begin{equation}\label{eq-3con}
\begin{split}
\phi_{0,\theta}&=F(\phi_0,\theta,z)+\lambda_0e^{\phi_1}\\
\phi_{1,\theta}&=F_r(\phi_0,\theta,z)-2\lambda_0 \phi_2e^{-\phi_1}+\lambda_1\\
\phi_{2,\theta}&=F_r(\phi_0,\theta,z)\phi_2+\frac12 F_{rr}(\phi_0,\theta,z)e^{2\phi_1}-2\lambda_1\phi_2+\lambda_2e^{\phi_1}.
\end{split}
\end{equation}
The first equation here is equivalent to the first equation in \eqref{eq-ODE}; 
the second has an additional term due to $\phi_2$; and, the third is new.
If we can find a branch of nontrivial solutions
of this system, then we can transform the general fast system
into the desired form up to and including ${\cal O}(\rho^2)$.

As above, we look for $2\pi$ periodic solutions. 
However, before extending the definition of ${\cal H}$,
we rewrite the third component of \eqref{eq-3con} as follows:
\begin{equation}\label{eq-3conp}
\frac{d}{d\theta}\left (\phi_{2}e^{-\phi_1}\right )
=\frac12 F_{rr}(\phi_0,\theta,z)e^{\phi_1}
+2\lambda_0\phi_2^2e^{-2\phi_1}-3\lambda_1\phi_2e^{-\phi_1}+\lambda_2.
\end{equation}
Now, we define ${\cal H}\; :\;\mathbb{R}\times\mathbb{R}^3\to \mathbb{R}^3$ 
in a manner similar to that employed in Step 1:
\begin{equation}\label{eq-defH3}
\begin{split}
&{\cal H}(z,\lambda_0,\lambda_1,\lambda_2)= \left (\begin{array}{c} {\cal H}_0(z,\lambda_0,\lambda_1,\lambda_2)\\{\cal H}_1(z,\lambda_0,\lambda_1,\lambda_2)\\{\cal H}_2(z,\lambda_0,\lambda_1,\lambda_2)\end{array}\right )=\\&
\left (\begin{array}{l}\int_0^{2\pi}F(\phi_0(\nu),\nu,z)d\nu+\lambda_0\int_0^{2\pi} e^{\phi_1(\nu)}d\nu\\
\int_0^{2\pi} (F_r(\phi_0(\nu),\nu,z)-2\lambda_0\phi_2(\nu)e^{-\phi_1(\nu)})d\nu +2\pi\lambda_1
\\ \int_0^{2\pi} (\frac12 F_{rr}(\phi_0(\nu),\nu,z)e^{\phi_1(\nu)}+2\lambda_0\phi_2(\nu)^2e^{-2\phi_1(\nu)}-3\lambda_1\phi_2(\nu)e^{-\phi_1(\nu)}+\lambda_2)d\nu\end{array}\right ).
\end{split}
\end{equation}

We analyse ${\cal H}$ in much the same manner as in Step 1.
Observe that, at $z=0$, we have $\lambda_0=0$ and $\lambda_1=0$. 
However, $\lambda_2$ is 
in general given by 
\[
\lambda_2=-\frac{1}{2\pi}\int_0^{2\pi} \frac12 F_{rr}(0,\nu,0)e^{\phi_1(\nu)}d\nu,
\]
and not zero. 

We now show that the off-diagonal elements 
of the Jacobian of ${\cal H}$ vanish.
First, an argument similar to that used in Step 1 shows that 
${\cal H}_{0,\lambda_1}(0,0,0,\lambda_{2,0})=0$. 
We also need to show that 
${\cal H}_{0,\lambda_2}(0,0,0,\lambda_{2,0})=0$
and ${\cal H}_{1,\lambda_2}(0,0,0,\lambda_{2,0})=0$.
To this end, 
we prove that $\phi_{0,\lambda_2}(0,0,0,\lambda_{2,0})\equiv 0$ and
$\phi_{1,\lambda_2}(0,0,0,\lambda_{2,0})\equiv 0$,
because then the identities 
${\cal H}_{0,\lambda_2}(0,0,0,\lambda_{2,0})=0$
and ${\cal H}_{1,\lambda_2}(0,0,0,\lambda_{2,0})=0$
follow in a straightforward way.
We carry out the proof for $\phi_{0,\lambda_2}$; the argument
for $\phi_{1,\lambda_2}$ is similar. 
Differentiating the first equation in \eqref{eq-3con} 
with respect to $\lambda_2$ 
and using the fact that $\lambda_0\; |_{z=0}=0$, we obtain
\begin{equation}\label{eq-keyeq}
\frac{d \phi_{0,\lambda_2}|_{z=0}}{d\theta}=F_r(0,\theta,0)\phi_{0,\lambda_2}|_{z=0}.
\end{equation}
The claim now follows from the assumption 
$\phi_0(0,z,\lambda_0,\lambda_1,\lambda_2)\equiv0$. 

Based on the above analysis,
it follows that
\[
\det(D{\cal H}(0,0,0,0))
={\cal H}_{0,\lambda_0}(0,0,0,0){\cal H}_{1,\lambda_1}(0,0,0,0)
        {\cal H}_{2,\lambda_2}(0,0,0,0),
\]
and an argument similar 
to the one used in Step 1
shows that this determinant is nonzero. 
In particular, 
${\cal H}_{0,\lambda_0}(0,0,0,0)$ and ${\cal H}_{1,\lambda_1}(0,0,0,0)$ 
are both $2\pi$ by a similar calculation.
To show that also ${\cal H}_{2,\lambda_2}(0,0,0,0)=2\pi$, 
we differentiate the third component in \eqref{eq-defH3} to obtain
\[
{\cal H}_{2,\lambda_2}(0,0,0,0)=\int_0^{2\pi} d\nu=2\pi.
\]
Hence, we may again use the Implicit Function Theorem
to conclude that there exists a branch of nontrivial solutions
of \eqref{eq-3con}, 
and the system may be put in the desired form
up to and including terms of $\rho^2$.
This completes the proof 
of Lemma~\ref{lem-chaco} for the case $N=2$.

\vskip0.3truein
\noindent
{\bf Step 3.}
In this third and final step of the proof,
we show that the lemma holds for general $N$
in \eqref{eq-chaco}.
We begin by writing
\begin{equation}\label{eq-chacogen}
F(r,\theta,z)=\sum_{j=0}^N \rho^j F_j(\phi_0,\phi_1,\ldots,\phi_N,\theta,z) +O(\rho^{N+1}),
\end{equation}
where $N>2$ is a natural number 
and $r$ and $\rho$ are related by formula \eqref{eq-chaco} 
associated to this choice of $N$.
For each $j$, 
the functions $F_j$ are complicated expressions 
involving $\phi_0$, $\phi_1$, $\ldots $, $\phi_j$. 
To simplify the notation,
we write $\Phi_j=(\phi_0,\phi_1,\ldots,\phi_j)$. We will give a more
precise description of the functions $F_j$ below. 
The equivalent of \eqref{eq-sysn2} is now
\begin{equation}\label{eq-sysnn}
\begin{split}
\dot\rho(1+\sum_{l=2}^Nl\rho^{l-1}\phi_le^{-\phi_1})=&(-\phi_{0,\theta}+F(\phi_0,\theta,z))e^{-\phi_1}+\rho(-\phi_{1,\theta}+ F_r(\phi_0,\theta,z))\\&+\sum_{j=2}^N \rho^j(-\phi_{j,\theta}+F_j(\Phi_j,\theta,z))e^{-\phi_1}.
\end{split}
\end{equation}
Let $\alpha_0=1$ and for each $j\ge 1$ set
\[
\alpha_j(\phi_0,\ldots ,\phi_n)=\frac{1}{j!}\partial^j_\rho \left (\frac{1}{1+\sum_{l=2}^nl\rho^{l-1}\phi_le^{-\phi_1}}\right )\big |_{\rho=0}.
\]
Note that, for given $j$, $\alpha_j$ depends on $\phi_0,\phi_1,\ldots,\phi_{j+1}$.
We also write
\begin{equation}\label{eq-E0E1dots}
\begin{split}
E_0&=(-\phi_{0,\theta}+F(\phi_0,\theta,z))e^{-\phi_1}\\
E_1&=-\phi_{1,\theta}+ F_r(\phi_0,\theta,z)\\
E_2&=(-\phi_{2,\theta}+F_2(\Phi_2,\theta,z))e^{-\phi_1}\\
&\vdots\\
E_N&=(-\phi_{N,\theta}+F_N(\Phi_N,\theta,z))e^{-\phi_1}.
\end{split}
\end{equation}
It follows that \eqref{eq-sysnn} 
maay be written in the following compact and insightful manner:
\begin{equation}\label{eq-sysnnn}
\dot\rho=\sum_{l=0}^N \rho^l \left(\sum_{j=0}^l \alpha_j E_{l-j}\right).
\end{equation}
We now define the set of equations
\begin{equation}\label{eq-equations1}
\begin{split}
E_0&=\lambda_0\\
E_1+\alpha_1E_0&=\lambda_1\\
E_2+\alpha_1E_1+\alpha_2E_0&=\lambda_2\\
&\vdots\\
\sum_{j=0}^N \alpha_j E_{N-j}&=\lambda_N.
\end{split}
\end{equation}
This enables us to rewrite 
\ref{eq-E0E1dots} as follows:
\begin{equation}\label{eq-equations2}
\begin{split}
E_0&=\lambda_0\\
E_1&=\lambda_1-\alpha_1\lambda_0\\
E_2&=\lambda_2-\alpha_1\lambda_1+(\alpha_1^2-\alpha_2)\lambda_0\\
&\vdots\\
E_N&=\sum_{j=0}^{N} \beta_j \lambda_{N-j},
\end{split}
\end{equation}
where $\beta_0=1$ and $\beta_1,\ldots,\beta_{N}$ are coefficients depending on $\alpha_0,\ldots,\alpha_N$.
Moreover, $\beta_j$ depends on $\alpha_0,\ldots,\alpha_j$ only.
Therefore, we have arrived at the 
following system of differential equations:
\begin{equation}\label{eq-ODEsn}
\begin{split}
\phi_{0,\theta}&=F(\phi_0,\theta,z)+\lambda_0e^{\phi_1}\\
\phi_{1,\theta}&=F_r(\phi_0,\theta,z)+\lambda_1-\alpha_1\lambda_0\\
\phi_{2,\theta}&=F_2(\Phi_2,\theta,z)
                  +(\lambda_2-\alpha_1\lambda_1
                     +(\alpha_1^2-\alpha_2)\lambda_0)e^{\phi_1}\\
&\vdots\\
\phi_{N,\theta}&=F_N(\Phi_N,\theta,z)
                 +\left(\sum_{j=0}^{N} \beta_j \lambda_{N-j}\right )e^{\phi_1},
\end{split}
\end{equation}
which is the analog for general $N$ of
the systems of differential equations
\eqref{eq-ODE} for $N=1$ 
and \eqref{eq-3con} for $N=2$.

Before defining ${\cal H}$,
we rewrite \eqref{eq-ODEsn} 
in a manner similar to that which was used above
to rewrite \eqref{eq-3con} 
(recall also \eqref{eq-3conp}).
Noting that 
\[
F_j(\Phi_j,\theta,z)=F_r(\phi_0,\theta,z)\phi_j+{\cal R}(\Phi_{j-1}),
\]
we replace the $j$th equation in \eqref{eq-ODEsn} by 
\begin{equation}\label{eq-jrepl}
\frac{d}{d\theta}\left (\phi_j e^{-\phi_1}\right )= {\cal R}(\Phi_{j-1})e^{-\phi_1}-(\lambda_1-\alpha_1\lambda_0)\phi_je^{-\phi_1}+\left(\sum_{l=0}^{j} \beta_l\lambda_{j-l}\right ).
\end{equation}
We will now define the function ${\cal H}$ 
in a manner analogous to \eqref{eq-HN1}
and \eqref{eq-defH3} used in STeps 1 and 2, respectively, for the cases
$N=1$ and $N=2$.
We let
\[
(\phi_0(\theta,z,\lambda_1,\ldots,\lambda_N),\phi_1(\theta,z,\lambda_1,\ldots,\lambda_N),\ldots,
\phi_N(\theta,z,\lambda_1,\ldots,\lambda_N))
\]
be the solutions of \eqref{eq-ODEsn} 
depending on the parameters $z$ and $\lambda_0,\ldots,\lambda_N$ 
that satisfy the initial conditions 
$\phi_j(0,z,\lambda_1,\ldots,\lambda_N)=0$, $j=0,1,\ldots, N$.
Further, we let ${\cal H}_0$, ${\cal H}_1$, and ${\cal H}_2$ 
be defined as in Step 2, and
let 
\begin{equation}\label{eq-defHj}
{\cal H}_j(z,\lambda_0,\lambda_1,\ldots,\lambda_N)=\int_0^{2\pi} \left({\cal R}(\Phi_{j-1})e^{-\phi_1}-(\lambda_1-\alpha_1\lambda_0)\phi_je^{-\phi_1}+\sum_{l=0}^{j} \beta_l \lambda_{j-l}\right ) d\nu,\qquad j=2,3,\ldots N.
\end{equation}
We first argue that we can solve the set of equations ${\cal H}_j(0,\lambda_0,\ldots,\lambda_N)$, $j=0,1,\ldots, N$, for a unique $N$-tuple $(\lambda_{0,0},\lambda_{1,0},\ldots,\lambda_{N,0})$.
Note that $\lambda_{0,0}=0$ and $\lambda_{1,0}=0$, 
by the same analysis used in Step 2. 
Hence,
\[
\sum_{k=0}^{j} \beta_k \lambda_{j-k}|_{z=0}=\sum_{k=0}^{j-2} \beta_k \lambda_{j-k}|_{z=0}.
\]
We argue by induction. Suppose that $\lambda_{0,0},\ldots, \lambda_{j-1,0}$ and the corresponding $\phi_0,\ldots, \phi_{j-1}$ are determined. 
Note that $\lambda_{j,0}$ must satisfy
\begin{equation}\label{eq-musa}
\lambda_{j,0}=-\int_0^{2\pi} \left({\cal R}(\Phi_{j-1})e^{-\phi_1}+\sum_{l=1}^{j-2} \beta_l \lambda_{j-l,0}\right ) d\nu.
\end{equation}
Since the right member of \eqref{eq-musa} depends only on $\lambda_0,\ldots,\lambda_{j-1}$ and $\phi_0,\ldots,\phi_{j-1}$ the value of $\lambda_{j,0}$ is uniquely determined.
Similarly, knowing $\lambda_{j,0}$, we can solve \eqref{eq-jrepl} for $\phi_j$ using \eqref{eq-jrepl} due to the fact that $\lambda_0= \lambda_1=0$, so that 
right member of \eqref{eq-jrepl}
depends only on  $\lambda_0,\ldots,\lambda_{j-1}$ 
and $\phi_0,\ldots,\phi_{j-1}$. 
Hence, $\phi_j$ is uniquely determined.

We now prove that 
$D_{\lambda_0,\ldots,\lambda_N}{\cal H}(0,\lambda_{0,0},\ldots,\lambda_{N,0})$ 
is non-singular. 
First, we prove that $\phi_{j,\lambda_k}|_{z=0}=0$ 
for any $j\in\{0,\ldots,N-1\}$ and $k\in\{j+1,\ldots,N\}$. 
The argument for $j=0$ and $j=1$ is analogous as in the case of $n=2$.
For general $j$, we proceed by induction, 
assuming that the claim holds for $0,1,\ldots, j-1$. 
The argument is, again, similar to that used in the case $N=2$.
We differentiate  \eqref{eq-jrepl} with respect to $\lambda_k$ 
and use the induction assumption,
the fact that $\beta_j$ is independent of $\phi_l$ for any $l>j+1$, and the fact that $\lambda_{0,0}=\lambda_{1,0}=0$. This gives
\begin{equation}\label{eq-jrepldiff}
\frac{d}{d\theta}\phi_{j,\lambda_k}(0,\lambda_{0,0},\ldots,\lambda_{N,0})= 0.
\end{equation}
By assumption, 
$\phi_j(0,z,\lambda_0,\ldots, \lambda_N)\equiv 0$. 
Hence, the claim follows. 
Now, differentiating \eqref{eq-defHj} and using a similar procedure, we obtain ${\cal H}_{j,\lambda_k}(0,\lambda_{0,0},\ldots,\lambda_{N,0})=0$ for $j\in\{0,\ldots,N-1\}$ and $k\in\{j+1,\ldots,N\}$.

It remains to prove that ${\cal H}_{j,\lambda_j}(0,\lambda_{0,0},\ldots,\lambda_{N,0})\neq 0$ for $j\in\{0,\ldots,N\}$. The proof for $j=0$ and $1$ is as above.
Let $j>1$.
Again by differentiating  \eqref{eq-defHj}, now with respect to
$\lambda_j$, and arguing analogously as above, we obtain
\begin{equation}\label{eq-jrepldiff-2}
{\cal H}_{j,\lambda_j}(0,\lambda_{0,0},\ldots,\lambda_{N,0})=\int_0^{2\pi}d\nu=2\pi.
\end{equation}
We can now apply 
the Implicit Function Theorem 
to obtain the functions $\lambda_0(z),\ldots,\lambda_j(z)$ as required.
This completes the third (and final) step of the proof of the lemma.

\hfill QED

\vspace{10pt}
\noindent {\bf Proof of Theorem \ref{th-red}} 
If we apply the sequence of transformations leading to 
\eqref{eq-aosysfps} to system \eqref{eq-osys},
we obtain a system of the form 
\begin{align}\label{eq-eaosysfps}
\begin{split}
\dot \rho&= G(\rho,z)+O(\rho^{N+1})+O(\eps)\\
\dot \theta& = 1\\
\dot z&=\eps \tilde g(\rho,\theta, z,\eps).
\end{split}
\end{align}
with 
\[
G(\rho,z)=\sum_{j=0}^{N-1}\lambda_j(z)\rho^j
\]
and the coefficient functions are as introduced in the proof of Lemma \ref{lem-chaco}. In particular hypotheses {\bf (H1)} and {\bf (H2)} imply $\lambda_0(0)=0$ and $\lambda_1(0)=0$. We now formulate the additional degeneracy assumptions in terms of the functions
$\lambda_j$. We will assume that $N\ge 3$.
\begin{equation}\label{eq-addhyp}
\begin{split}
\frac{d\lambda_0}{dz}\left (0\right )&\neq 0.\\
\lambda_2(0)&\neq 0.
\end{split}
\end{equation}
If \eqref{eq-addhyp} holds in addition to {\bf (H1)} and {\bf (H2)} we replace the variable $z$ by $\tilde z=\lambda_0(z)$ and perform scalings and translations to arrive at \eqref{eq-vdpl}, taking $\eps$ sufficiently small as necessary.\hfill QED

\vspace{10pt}
We now make some remarks relating the normal form equation \eqref{eq-vdpl} derived in this section and the forced van der Pol  \eqref{fvdp}. Clearly, the fast subsystems are similar. If in addition to \eqref{eq-addhyp} we assume that $\lambda_3(0)\neq 0$ and make some
assumptions about the signs of the coefficients, then the two fast systems are the same to lowest order (up to a simple transformation). The situation with the slow equation is more complicated. In systems with time-periodic forcing, we have shown that \eqref{fvdp} is a normal form. Moreover, because forced systems are not generic in the larger class of general slow/fast dynamical systems, we expect the dynamics in this larger class to be even richer. In addition, if we assume that $\tilde g(0,\theta_0,0,0)=0$ for some $\theta_0\in [0,\, 2\pi)$
then we may obtain some canard dynamics. In general these slow/fast systems will have folded singularities and associated canard dynamics. %This case is similar to   \eqref{fvdp} with $ a\le |b| $ and  is not studied in our paper. However the assumption $\tilde g(0,\theta, 0)\equiv 0$, equivalent to the hypothesis $b=0$ in \eqref{fvdp} which is  assumed in our paper,  
%is not generic in the context of  \eqref{eq-vdpl} within any class of vector fields of finite codimension $k$, for any integer $k$. 
%Hence, it appears that the case we study in this paper is more restrictive than the general case of torus canards.

%---------------------------------------------------------------------------------
\section{Conclusions and Discussion}		\label{sec:discussion}
%---------------------------------------------------------------------------------

In this article, we have established the existence of a number of different types of canard solutions of the forced van der Pol equation
\eqref{fvdp} across the entire range of forcing frequencies $\omega>0$. 
Most interestingly, we have found numerically that the families of primary maximal canards and maximal torus canards are organised along single branches in parameter space. 
In the low-frequency regime ($\omega={\cal O}(\eps)$), Theorem \ref{thm:lowfrequency} demonstrates the existence of the primary maximal canards of the FSN I points and establishes that formula \eqref{mainformulas-Reg1} gives the loci of the folds of the primary maximal canards in the $(a,b,\overline{\omega})$ parameter space with $b={\cal O}(\sqrt{\eps})$. In the intermediate-frequency ($\omega={\cal O}(\sqrt{\eps})$) and high-frequency regime ($\omega={\cal O}(1)$), Theorem~\ref{thm:intermediatefrequency} and Corollary~\ref{cor:highfrequency}, respectively, establish the existence of maximal torus canards as well as the formulas \eqref{mainformulas-Reg2} and \eqref{mainformulas-Reg3}, which explicitly give the locations of the folds of the maximal torus canards in the $(a,b,\omega)$ parameter space with $b={\cal O}(\eps)$. These maximal torus canards lie precisely in the intersection of the persistent critical manifolds of attracting limit cycles and of repelling limit cycles. They are the analogs in one-higher-dimension of the maximal headless ducks of the unforced van der Pol equation, see for example \cite{Benoit1981,Dumortier1996,Eckhaus1983}. Moreover, they are similar to the folds of maximal torus canards observed earlier in a rotated system of ven der Pol type, see Figure 5 in \cite{Benes2011}.

It was also shown that these analytical results are all representations of the same formula \eqref{mainformulas} that holds across the entire range of forcing frequencies $\omega>0$ for the appropriate values of $b$, and that these formulas agree well with the results obtained from numerical continuations over the parameter regions in which they apply. Moreover, in the limit $\omega\to
\infty$, the torus canards appear to be rotated copies of the limit cycle canards
that exist in the planar unforced van der Pol equation, and the interval of
$a$ values for which the maximal torus canards exist shrinks to the
value $a_c(\eps)=1-\frac{\eps}{8}$ at which the maximal headless duck
solution exists in the unforced equation, recall \cite{Baer1986,
Benoit1981, Braaksma1993, Dumortier1996, Eckhaus1983}.

It is worth noting that the analytical results presented here for
the torus canards of \eqref{fvdp} agree with and expand upon by the general topological
analysis presented in Section 6 of \cite{Burke2012}. There, fast-slow
systems with two fast variables and one slow variable were studied in
which there is a torus canard explosion in the transition regime from
stable periodic spiking (tonic spiking) to bursting, and examples were
given, including of the Hindmarsh-Rose equation, the Morris-Lecar-Terman
model, and the Wilson-Cowan-Izhikevich system. In particular, it was shown
using topological arguments that there must be a sequence of torus canards
in these transition regimes in order to satisfy the property of continuous
dependence of solutions on parameters. The topological analysis presented
in \cite{Burke2012} is the analog in one higher dimension of the
topological analysis first used in \cite{Benoit1981} to establish the
existence of an explosion of limit cycle canards in the transition between
asymptotically stable solutions and full-blown relaxation oscillations in
the unforced, planar van der Pol equation.

In this article, we also studied the branches of the secondary maximal canards, which exist
in the low-frequency regime in \eqref{fvdp}. Secondary canards lie close
to the primary strong canard for most of their lengths, and in addition
they make finitely many loops near the bottom of $\Gamma$, recall Figure
7. They are indexed by the number of loops and by the height in the $y$
variable of the jumps from the repelling slow manifold back to the
attracting slow manifold. We showed how the dynamics of these secondary
canards changes as the parameters change along the fold curves, and we
identified the mechanism by which these branches turn around well before
they get into the high-frequency regime. In particular, the turning points
correspond precisely to the parameter values at which the fold curve of
maximal canards is crossed and two intersection points of the attracting
and repelling slow manifolds are created (or annihilated). In addition, we
identified how new LAO segments are added to the secondary canard
solutions at points at which the direction of spiralling of the repelling
slow manifold is reversed, recall Figure 9.

Finally, we proved that the fvdP equation \eqref{fvdp} is a normal form for a class of slow/fast systems with two fast variables and one slow variable, which possess a non-degenerate fold of limit cycles in the fast system and which exhibit the torus canard explosion phenomenon. Thus, the methods and results obtained here for \eqref{fvdp} extend naturally to a large class of slow/fast systems with single-frequency time-periodic forcing.

To conclude this article, we discuss a number of topics related to the
canard solutions of the fvdP \eqref{fvdp}. First, the fold curves of the primary maximal canards in the low-frequency regime and the fold
curves of the maximal torus canards in the intermediate- and high-frequency regimes together serve as the boundary of the MMO regime in \eqref{fvdp}, recall Figures~\ref{fig:motivation} and \ref{fig:mmos}.

Second, there are many different branches of resonance curves, curves of torus bifurcations, saddle-nodes of periodic orbits, period-doubling curves of periodic orbits, and so forth, all of which lie in the interior of the MMO region in parameter space. These bifurcation curves play important roles as the boundaries between orbit segments with different numbers of SAOs and LAOs. These are the subject of ongoing research.

%\newpage
%------------------------------------------------------------
\section*{Acknowledgments} 
%------------------------------------------------------------
The research of J.B., T.J.K., and T.V. was partially supported by NSF-DMS 1109587. T.J.K. thanks INRIA for their hospitality and for providing a climate conducive to research and collaboration during a visit. The authors thank Mark Kramer, Nick Benes, John Mitry and Martin Wechselberger for useful conversations.

%------------------------------------------------------------------------------------------------------------------------
\appendix
%------------------------------------------------------------------------------------------------------------------------
%---------------------------------------------------------------------------------
\section{Proof of the Existence of a Torus Bifurcation} 	\label{sec:omegafast}
%---------------------------------------------------------------------------------
In this appendix, we consider the forced vdP oscillator \eqref{fvdp} subject to high-frequency forcing ($\omega = \mathcal{O}(1)$). In this case, $y$ is the only slow variable, and there is no critical manifold. As such, the fast dynamics are dominant throughout the phase space, and system \eqref{fvdp} can be interpreted as a regularly perturbed, non-autonomous problem. 

We establish the existence of a torus bifurcation using second-order averaging and equivariant normal form theory in Section \ref{subsec:omegafastav}, and we present the asymptotics of the torus bifurcation parameter value $a_{\rm TB}(b,\omega,\varepsilon)$ in Section \ref{subsec:omegafastasymp}.

%------------------------------------------
\subsection{Second-Order Averaging and a Torus Bifurcation in the Regime $\omega = \mathcal{O}(1)$}	\label{subsec:omegafastav}
%------------------------------------------

In this section, we perform a standard second-order averaging analysis of system \eqref{fvdp} and use equivariant normal form theory in the regime $\omega=\mathcal{O}(1)$ to demonstrate that there exists a smooth function $a=a_{\rm TB}(b,\omega,\varepsilon)$ at which the system possesses a torus bifurcation
in which a stable two-torus is born. 
First, we change variables so that the fold is located at the origin: $({\tilde x},{\tilde y}) = (1-x,y+\frac{2}{3})$
and $\alpha = a - 1$. The forced vdP equation transforms to 
\begin{equation*}
\begin{aligned}
\tilde{x}^\prime &=-{\tilde y}+{\tilde x}^2-\frac{1}{3}{\tilde x}^3\\
{\tilde y}^\prime &=\varepsilon\left( {\tilde x}+\alpha+b\cos\theta \right)\\
\theta^\prime &=\omega.
\end{aligned}
\end{equation*}
Then, to carry out the second-order averaging analysis, it is natural to use the following scaling, which comes from the central chart
of the desingularization, or blow-up, analysis 
used in Section \ref{sec:AnalysisRegion2+3} and in~\cite{Krupa2001}:
\[ {\tilde x}=\sqrt{\eps}\, \overline{x}, \quad {\tilde y}=\eps\overline{y}, \quad b=\sqrt{\eps}\, \overline{b}, \quad \alpha=\sqrt{\eps} \, \overline{\alpha}, \]
and to rescale time by $t \longrightarrow \omega t$. 
Hence, after dropping the bars, 
the system has the following non-autonomous form:
\begin{equation}        \label{eq:rescaled_non-autonomous}
\begin{split}
\dot{x}&=\frac{\delta}{\omega}\left(-y+x^2 \right)-\frac{\delta^2}{3\omega}x^3\\
\dot{y}&=\frac{\delta}{\omega}\left(x+\alpha+b\cos (t_0+t) \right),
\end{split}
\end{equation}
where $\delta=\sqrt{\eps}$. The choice of $t_0$ has no effect on the analysis. 

Next, we apply the near-identity change of variables used in second-order averaging \cite{SV1982}, so that system~\eqref{eq:rescaled_non-autonomous} transforms into
\begin{equation}        \label{eq:rescaled_non-autonomous_acv}
\begin{split}
\dot{\xi_1}&=\frac{\delta}{\omega}\left(-\xi_2+\xi_1^2 \right) -\frac{\delta^2}{3\omega}\xi_1^3+\tilde{G}(\xi_1,\xi_2,t,\delta)\\
\dot{\xi_2}&=\frac{\delta}{\omega}\left(\xi_1+\alpha\right)+O(\delta^3),
\end{split}
\end{equation}
with $\tilde{G}(\xi_1,\xi_2,t,\delta)=O(\delta^3)$. We label \eqref{eq:rescaled_non-autonomous_acv} as the `intermediate' system; it is smoothly conjugate to the original system.

We are interested in the averaged system,
\begin{equation}        \label{eq:averaged_system}
\begin{split}
\dot{\bar{x}}&=\frac{\delta}{\omega}\left(-\bar{y}+\bar{x}^2 \right)-\frac{\delta^2}{3\omega}\bar{x}^3\\
\dot{\bar{y}}&=\frac{\delta}{\omega}\left(\bar{x}+\alpha\right).
\end{split}
\end{equation}
This system has a unique $S^1$ equivariant normal form \cite{GSS1988} due to its symmetry properties. Moreover, the time-$T$ map of this normal form must be the $S^1$ normal form of the time-$T$ map, since the two operations commute and since $S^1$ equivariant normal forms
are unique. Let $\ts$ denote the Poincare map of this normal form. At $\alpha=0$, the eigenvalues of the map $\ts$ satisfy the non-resonance condition: they are not equal to the first four strong resonant eigenvalues. 
Also, the second Liapunov coefficient is negative; 
in fact, for the averaged system, 
the second Liapunov coefficient is known to be $K\delta$, 
where $K<0$ is a constant. 
Hence, at $\alpha=0$, the map $\ts$
satisfies the basic hypotheses of the Hopf bifurcation for maps; see conditions A) and B), respectively, and Theorem 2 of~\cite{Lanford1973}.
Therefore, at $\alpha=0$, 
the map $\ts$ undergoes a non-degenerate, 
supercritical Hopf bifurcation in which a normally-hyperbolic invariant circle 
is created. 
Hence, one also sees directly that the averaged system 
\eqref{eq:averaged_system} undergoes a torus bifurcation at $\alpha=0$ 
in which the limit cycle 
becomes unstable and a stable invariant torus is created.

We now demonstrate that the full system 
\eqref{eq:rescaled_non-autonomous} 
undergoes a torus bifurcation 
at some $\alpha_{\rm TB}$ near $\alpha=0$, 
in which a stable invariant torus is created. 
In particular, we show that the Poincare map $\s$ 
of the full system \eqref{eq:rescaled_non-autonomous} 
possesses an invariant circle $\delta$-close to the one of $\ts$. 

First, we observe that the same near-identity coordinate change employed above to put the averaged system \eqref{eq:averaged_system}
into its $S^1$ equivariant normal form also puts the intermediate system \eqref{eq:rescaled_non-autonomous_acv} into its $S^1$ equivariant normal form, up to and including ${\cal O}(\delta^2)$. Let $\s_i(z)$ denote the Poincare map of this normal form. Again, this map must be the time-$T$ map of the normal form, due to uniqueness in the $S^1$ equivariant case. Next, we observe that, by standard second-order averaging theory, the Poincare maps are close, i.e., on time intervals of length $\mathcal{O}(1/\delta)$, 
\[ \left| \s_i(z)-\ts(z) \right| = {\cal O}(\delta). \]
In fact, for each $\alpha$ sufficiently small, the map $\s_i$ is part of exactly the type of one-parameter family of maps studied in \cite{Lanford1973}, with the map of the averaged system being the `unperturbed' map. 
Hence, for some $\alpha$ near $\alpha=0$, the map $\s_i$ also undergoes a non-degenerate, super-critical Hopf bifurcation in which a normally-attracting invariant circle is created.

Finally,
we observe that
the time-$T$ maps $\s$ and $\s_i$
are smoothly conjugate.
Hence, it follows that
the map $\s$ of the original
system~\eqref{eq:rescaled_non-autonomous}
also has a non-degenerate Hopf bifurcation,
and therefore that
the original vector field
\eqref{eq:rescaled_non-autonomous}
has a torus bifurcation
at some $\alpha_{\rm TB}$ near zero,
in which an attracting invariant torus
is created.
This concludes the demonstration.

\begin{remark}
One boundary of the parameter $\alpha$
for the existence
of an invariant torus
for the full system
is given by the birth of the canard regime.
\end{remark}

%------------------------------------------
\subsection{Asymptotic Expansion of $a_{\rm TB}(b, \omega, \varepsilon)$}	\label{subsec:omegafastasymp}
%------------------------------------------

In this section, we present the asymptotic expansion of 
$a = a_{TB}(b,\omega,\eps)$ for small $b$. 
The unforced vdP equation
has an equilibrium at $(x,y)=(a,f(a))$,
which undergoes a Hopf bifurcation at $a=1$.
This corresponds
to a periodic orbit of \eqref{fvdp}
at $b=0$
which undergoes a torus bifurcation at $a=1$.
We seek periodic solutions
of \eqref{fvdp} as an asymptotic series in $b$,
i.e., let
\begin{align*}
x(t) = \sum_{k=0}^{\infty} b^k x_k(t), \quad
y(t) = \sum_{k=0}^{\infty} b^k y_k(t).
\end{align*}
Substitution into \eqref{fvdp} yields
\begin{align*}
\dot{x}_0 &= y_0 - f(x_0), \\
\dot{y}_0 &= \eps (a-x_0),
\end{align*}
at leading order,
which is the planar van der Pol equation.
This has an equilibrium
at $(x,y)=(a,f(a))$.
The $\mathcal{O}(b^1)$ system is
\begin{align*}
\dot{x}_1 &= y_1 - f^\prime(x_0) \,x_1, \\
\dot{y}_1 &= \eps (-x_1+\cos \omega t).
\end{align*}
The $\mathcal{O}(b^1)$ system is linear with solutions that are linear combinations of $\cos \omega t, \sin \omega t$, and $\exp (\lambda t)$,
where $2\lambda = 1-a^2 \pm \sqrt{(1-a^2)^2-4\eps}$.

\begin{remark}
Note that when $a=1$ and $\omega = \sqrt{\eps}$, there is a resonance.
\end{remark}

When a damped oscillator is driven with a periodic forcing function, the result may be a periodic response at the same frequency as the forcing function (see Figure \ref{fig:mmos}(a)). Since the unforced oscillation is dissipated due to the damping, it is absent from the steady state behaviour. Thus, we seek periodic solutions of the $\mathcal{O}(b^1)$ system of period $T = \frac{2\pi}{\omega}$. The solutions are
\begin{align*}
x_1(t) &= \frac{\left(a^2-1\right) \varepsilon  \omega  \sin (t \omega )+\varepsilon  \left(\varepsilon -\omega ^2\right) \cos (t \omega )}{\left(a^2-1\right)^2 \omega^2+\left(\varepsilon -\omega ^2\right)^2}, \\
y_1(t) &= \frac{\varepsilon  \left(\left(a^2-1\right) \varepsilon  \cos (t \omega )+\omega  \left(a^4-2 a^2-\varepsilon +\omega ^2+1\right) \sin (t \omega)\right)}{\left(a^2-1\right)^2 \omega ^2+\left(\varepsilon -\omega ^2\right)^2}.
\end{align*}

We now compute the stability of a periodic solution of \eqref{fvdp}. Let $(x_\gamma,y_\gamma)$ denote the periodic solution and let $(u,v)=(x-x_\gamma,y-y_\gamma)$ be a perturbation of this orbit. Then the perturbations evolve according to
\[ \begin{pmatrix} \dot{u} \\ \dot{v} \end{pmatrix} = Df(x_\gamma,y_\gamma) \begin{pmatrix} u \\ v \end{pmatrix} - \begin{pmatrix} \frac{1}{2}f^{\prime \prime}(x_\gamma)u^2 + \frac{1}{6} f^{\prime \prime \prime} (x_\gamma) u^3 \\ 0 \end{pmatrix}, \]
where the Jacobian evaluated along $(x_\gamma,y_\gamma)$ is the $T$-periodic matrix:
\[ Df(x_\gamma,y_\gamma) = \begin{pmatrix} -f^\prime(x_\gamma) & 1 \\ -\eps & 0 \end{pmatrix}. \]
The Floquet multipliers $\rho_1$ and $\rho_2$ satisfy
\[ \rho_1 \rho_2 = \exp \left( \int_0^T \operatorname{tr} Df(x_\gamma,y_\gamma)\, dt \right). \]
A Neimark-Sacker bifurcation occurs when the multipliers are $\rho = e^{\pm i \mu}$ for some $\mu$. That is, a torus bifurcation occurs when $\int_0^T \operatorname{tr} Df(x_\gamma,y_\gamma)\, dt = 0$. This gives the following relation between $a,b,\eps$, and $\omega$ for the location of the torus bifurcation:
\begin{equation}        \label{TB}
1-a^2-\frac{1}{2} \frac{b^2 \eps^2}{(a^2-1)^2 \omega^2+(\eps-\omega^2)^2} =0.
\end{equation}
Comparisons between the theoretical prediction above and the results of numerical continuation simulations show good agreement for the indicated parameter regions. Divergence from the above perturbation analysis occurs precisely at the resonances (figures not shown).

\begin{remark}
Note that \eqref{TB} is accurate up to terms of order $\mathcal{O}(b^3)$. That is, there are no order $\mathcal{O}(b^2)$ corrections in \eqref{TB} from the $b^2 x_2$ terms in the asymptotic expansion due to symmetry. More precisely, symmetry considerations give
\[ -2ab^2 \int_0^T x_2 \, dt = 0. \]
\end{remark}

%%%%%%%%      REFERENCES      %%%%%%%%

%=========================================================================================
\end{document}